\magnification\magstep1
\baselineskip15pt
\newread\AUX\immediate\openin\AUX=\jobname.aux
\def\ref#1{\expandafter\edef\csname#1\endcsname}
\ifeof\AUX\immediate\write16{\jobname.aux gibt es nicht!}\else
\input \jobname.aux
\fi\immediate\closein\AUX
\def\today{\number\day.~\ifcase\month\or
  Januar\or Februar\or M{\"a}rz\or April\or Mai\or Juni\or
  Juli\or August\or September\or Oktober\or November\or Dezember\fi
  \space\number\year}
\font\sevenex=cmex7
\scriptfont3=\sevenex
\font\fiveex=cmex10 scaled 500
\scriptscriptfont3=\fiveex
\def\A{{\bf A}}

\def\P{{\bf P}}

\def\XS{{\widetilde X}}

\def\phi{\varphi}
\def\epsilon{\varepsilon}
\def\theta{\vartheta}
\def\uauf{\lower1.7pt\hbox to 3pt{%
\vbox{\offinterlineskip
\hbox{\vbox to 8.5pt{\leaders\vrule width0.2pt\vfill}%
\kern-.3pt\hbox{\lams\char"76}\kern-0.3pt%
$\raise1pt\hbox{\lams\char"76}$}}\hfil}}
\def\cite#1{\expandafter\ifx\csname#1\endcsname\relax
{\bf?}\immediate\write16{#1 ist nicht definiert!}\else\csname#1\endcsname\fi}
\def\expandwrite#1#2{\edef\next{\write#1{#2}}\next}
\def\neverexpand{\noexpand\noexpand\noexpand}
\def\strip#1\ {}
\def\ncite#1{\expandafter\ifx\csname#1\endcsname\relax
{\bf?}\immediate\write16{#1 ist nicht definiert!}\else
\expandafter\expandafter\expandafter\strip\csname#1\endcsname\fi}
\newwrite\AUX
\immediate\openout\AUX=\jobname.aux
\newcount\Abschnitt\Abschnitt0
\def\beginsection#1. #2 \par{\advance\Abschnitt1%
\vskip0pt plus.10\vsize\penalty-250
\vskip0pt plus-.10\vsize\bigskip\vskip\parskip
\edef\TEST{\number\Abschnitt}
\expandafter\ifx\csname#1\endcsname\TEST\relax\else
\immediate\write16{#1 hat sich geaendert!}\fi
\expandwrite\AUX{\neverexpand\ref{#1}{\TEST}}
\leftline{\bf\number\Abschnitt. \ignorespaces#2}%
\nobreak\smallskip\noindent\SATZ1}
\def\Proof:{\par\noindent{\it Proof:}}
\def\Remark:{\ifdim\lastskip<\medskipamount\removelastskip\medskip\fi
\noindent{\bf Remark:}}
\def\Remarks:{\ifdim\lastskip<\medskipamount\removelastskip\medskip\fi
\noindent{\bf Remarks:}}
\def\Definition:{\ifdim\lastskip<\medskipamount\removelastskip\medskip\fi
\noindent{\bf Definition:}}
\def\Example:{\ifdim\lastskip<\medskipamount\removelastskip\medskip\fi
\noindent{\bf Example:}}
\newcount\SATZ\SATZ1
\def\proclaim #1. #2\par{\ifdim\lastskip<\medskipamount\removelastskip
\medskip\fi
\noindent{\bf#1.\ }{\it#2}\Par
\ifdim\lastskip<\medskipamount\removelastskip\goodbreak\medskip\fi}
\def\Aussage#1{%
\expandafter\def\csname#1\endcsname##1.{\ifx?##1?\relax\else
\edef\TEST{#1\penalty10000\ \number\Abschnitt.\number\SATZ}
\expandafter\ifx\csname##1\endcsname\TEST\relax\else
\immediate\write16{##1 hat sich geaendert!}\fi
\expandwrite\AUX{\neverexpand\ref{##1}{\TEST}}\fi
\proclaim {\number\Abschnitt.\number\SATZ. #1\global\advance\SATZ1}.}}
\Aussage{Theorem}
\Aussage{Proposition}
\Aussage{Corollary}
\Aussage{Lemma}
\font\la=lasy10
\def\strich{\hbox{$\vcenter{\hbox
to 1pt{\leaders\hrule height -0,2pt depth 0,6pt\hfil}}$}}
\def\dashedrightarrow{\hbox{%
\hbox to 0,5cm{\leaders\hbox to 2pt{\hfil\strich\hfil}\hfil}%
\kern-2pt\hbox{\la\char\string"29}}}

\def\Bindestrich{\penalty10000-\hskip0pt}
\let\_=\Bindestrich
\def\.{{\sfcode`.=1000.}}
\def\Links#1{\llap{$\scriptstyle#1$}}
\def\Rechts#1{\rlap{$\scriptstyle#1$}}
\def\Par{\par}
\def\:={\mathrel{\raise0,9pt\hbox{.}\kern-2,77779pt
\raise3pt\hbox{.}\kern-2,5pt=}}
\def\=:{\mathrel{=\kern-2,5pt\raise0,9pt\hbox{.}\kern-2,77779pt
\raise3pt\hbox{.}}} \def\mod{/\mskip-5mu/}
\def\into{\hookrightarrow}
\def\pfeil{\rightarrow}
\def\untenPf{\downarrow}

\def\pf#1{\buildrel#1\over\rightarrow}

\def\Ugleich{\hbox{$\cup$\kern.5pt\vrule depth -0.5pt}}
\def\|#1|{\mathop{\rm#1}\nolimits}
\def\<{\langle}
\def\>{\rangle}
\let\Times=\times
\def\times{\mathop{\Times}}
\let\Otimes=\otimes
\def\otimes{\mathop{\Otimes}}
\catcode`\@=11
\def\hex#1{\ifcase#1 0\or1\or2\or3\or4\or5\or6\or7\or8\or9\or A\or B\or
C\or D\or E\or F\else\message{Warnung: Setze hex#1=0}0\fi}
\def\fontdef#1:#2,#3,#4.{%
\alloc@8\fam\chardef\sixt@@n\FAM
\ifx!#2!\else\expandafter\font\csname text#1\endcsname=#2
\textfont\the\FAM=\csname text#1\endcsname\fi
\ifx!#3!\else\expandafter\font\csname script#1\endcsname=#3
\scriptfont\the\FAM=\csname script#1\endcsname\fi
\ifx!#4!\else\expandafter\font\csname scriptscript#1\endcsname=#4
\scriptscriptfont\the\FAM=\csname scriptscript#1\endcsname\fi
\expandafter\edef\csname #1\endcsname{\fam\the\FAM\csname text#1\endcsname}
\expandafter\edef\csname hex#1fam\endcsname{\hex\FAM}}
\catcode`\@=12 

\fontdef Ss:cmss10,,.
\fontdef Fr:eufm10,eufm7,eufm5.
\def\fa{{\Fr a}}
\def\fb{{\Fr b}}
\def\fc{{\Fr c}}

\def\fg{{\Fr g}}
\def\fh{{\Fr h}}

\def\fk{{\Fr k}}
\def\fl{{\Fr l}}
\def\fm{{\Fr m}}

\def\fp{{\Fr p}}
\def\fq{{\Fr q}}\def\fQ{{\Fr Q}}
\def\fr{{\Fr r}}

\def\ft{{\Fr t}}
\def\fu{{\Fr u}}

\fontdef bbb:msbm10,msbm7,msbm5.
\fontdef mbf:cmmib10,cmmib7,.
\fontdef msa:msam10,msam7,msam5.
\def\CC{{\bbb C}}

\def\NN{{\bbb N}}
\def\RR{{\bbb R}}

\def\cA{{\cal A}}\def\cB{{\cal B}}\def\cC{{\cal C}}
\def\cE{{\cal E}}\def\cH{{\cal H}}

\def\cO{{\cal O}}
\def\cR{{\cal R}}
\def\cU{{\cal U}}

\mathchardef\leer=\string"0\hexbbbfam3F
\mathchardef\subsetneq=\string"3\hexbbbfam24
\mathchardef\semidir=\string"2\hexbbbfam6E
\mathchardef\dirsemi=\string"2\hexbbbfam6F
\mathchardef\haken=\string"2\hexmsafam78
\mathchardef\auf=\string"3\hexmsafam10
\let\OL=\overline
\def\overline#1{{\hskip1pt\OL{\hskip-1pt#1\hskip-1pt}\hskip1pt}}

\def\Bq{{\overline{B}}}
\def\Cq{{\overline{C}}}
\def\Dq{{\overline{D}}}

\def\fQ{{\overline{f}}}

\def\hq{{\overline{h}}}

\def\Kq{{\overline{K}}}

\def\Mq{{\overline{M}}}

\def\Sq{{\overline{S}}}
\def\Tq{{\overline{T}}}
\def\Uq{{\overline{U}}}
\def\Vq{{\overline{V}}}

\def\Xq{{\overline{X}}}

\def\zq{{\overline{z}}}
%
\abovedisplayskip 9.0pt plus 3.0pt minus 3.0pt
\belowdisplayskip 9.0pt plus 3.0pt minus 3.0pt
\newdimen\Grenze\Grenze2\parindent\advance\Grenze1em
\newdimen\Breite
\newbox\DpBox
\def\NewDisplay#1$${\Breite\hsize\advance\Breite-\hangindent
\setbox\DpBox=\hbox{\hskip2\parindent$\displaystyle{#1}$}%
\ifnum\predisplaysize<\Grenze\abovedisplayskip\abovedisplayshortskip
\belowdisplayskip\belowdisplayshortskip\fi
\global\futurelet\nexttok\WEITER}
\def\WEITER{\ifx\nexttok\qed\expandafter\leftQEDdisplay
\else\leftdisplay\fi}
\def\leftdisplay{\hskip-\hangindent\leftline{\box\DpBox}$$}
\def\leftQEDdisplay{\hskip-\hangindent
\line{\copy\DpBox\hfill\lower\dp\DpBox\copy\QEDbox}%
\belowdisplayskip0pt$$\bigskip\let\nexttok=}
\everydisplay{\NewDisplay}
\newbox\QEDbox
\newbox\nichts\setbox\nichts=\vbox{}\wd\nichts=2mm\ht\nichts=2mm
\setbox\QEDbox=\hbox{\vrule\vbox{\hrule\copy\nichts\hrule}\vrule}
\def\qed{\leavevmode\unskip\hfil\null\nobreak\hfill\copy\QEDbox\medbreak}
\newdimen\HIindent
\newbox\HIbox
\def\setHI#1{\setbox\HIbox=\hbox{#1}\HIindent=\wd\HIbox}
\def\HI#1{\par\hangindent\HIindent\hangafter=0\noindent\leavevmode
\llap{\hbox to\HIindent{#1\hfil}}\ignorespaces}

\newdimen\maxSpalbr
\newdimen\altSpalbr

\def\beginrefs{%
\expandafter\ifx\csname Spaltenbreite\endcsname\relax
\def\Spaltenbreite{1cm}\immediate\write16{Spaltenbreite undefiniert!}\fi
\expandafter\altSpalbr\Spaltenbreite
\maxSpalbr0pt
\def\L|Abk:##1|Sig:##2|Au:##3|Tit:##4|Zs:##5|Bd:##6|S:##7|J:##8||{%
\edef\TEST{[##2]}
\expandafter\ifx\csname##1\endcsname\TEST\relax\else
\immediate\write16{##1 hat sich geaendert!}\fi
\expandwrite\AUX{\neverexpand\ref{##1}{\TEST}}
\setHI{[##2]\ }
\ifnum\HIindent>\maxSpalbr\maxSpalbr\HIindent\fi
\ifnum\HIindent<\altSpalbr\HIindent\altSpalbr\fi
\HI{[##2]}
\ifx-##3\relax\else{##3}: \fi
\ifx-##4\relax\else{##4}{\sfcode`.=3000.} \fi
\ifx-##5\relax\else{\it ##5\/} \fi
\ifx-##6\relax\else{\bf ##6} \fi
\ifx-##8\relax\else({##8})\fi
\ifx-##7\relax\else, {##7}\fi\Par}
\def\B|Abk:##1|Sig:##2|Au:##3|Tit:##4|Reihe:##5|Verlag:##6|Ort:##7|J:##8||{%
\edef\TEST{[##2]}
\expandafter\ifx\csname##1\endcsname\TEST\relax\else
\immediate\write16{##1 hat sich geaendert!}\fi
\expandwrite\AUX{\neverexpand\ref{##1}{\TEST}}
\setHI{[##2]\ }
\ifnum\HIindent>\maxSpalbr\maxSpalbr\HIindent\fi
\ifnum\HIindent<\altSpalbr\HIindent\altSpalbr\fi
\HI{[##2]}
\ifx-##3\relax\else{##3}: \fi
\ifx-##4\relax\else{##4}{\sfcode`.=3000.} \fi
\ifx-##5\relax\else{(##5)} \fi
\ifx-##7\relax\else{##7:} \fi
\ifx-##6\relax\else{##6}\fi
\ifx-##8\relax\else{ ##8}\fi\Par}
\def\Pr|Abk:##1|Sig:##2|Au:##3|Artikel:##4|Titel:##5|Hgr:##6|Reihe:{%
\edef\TEST{[##2]}
\expandafter\ifx\csname##1\endcsname\TEST\relax\else
\immediate\write16{##1 hat sich geaendert!}\fi
\expandwrite\AUX{\neverexpand\ref{##1}{\TEST}}
\setHI{[##2]\ }
\ifnum\HIindent>\maxSpalbr\maxSpalbr\HIindent\fi
\ifnum\HIindent<\altSpalbr\HIindent\altSpalbr\fi
\HI{[##2]}
\ifx-##3\relax\else{##3}: \fi
\ifx-##4\relax\else{##4}{\sfcode`.=3000.} \fi
\ifx-##5\relax\else{In: \it ##5}. \fi
\ifx-##6\relax\else{(##6)} \fi\PrII}
\def\PrII##1|Bd:##2|Verlag:##3|Ort:##4|S:##5|J:##6||{%
\ifx-##1\relax\else{##1} \fi
\ifx-##2\relax\else{\bf ##2}, \fi
\ifx-##4\relax\else{##4:} \fi
\ifx-##3\relax\else{##3} \fi
\ifx-##6\relax\else{##6}\fi
\ifx-##5\relax\else{, ##5}\fi\Par}
\bgroup
\baselineskip12pt
\parskip2.5pt plus 1pt
\hyphenation{Hei-del-berg}
\sfcode`.=1000
\beginsection References. References

}
\def\endrefs{%
\expandwrite\AUX{\neverexpand\ref{Spaltenbreite}{\the\maxSpalbr}}
\ifnum\maxSpalbr=\altSpalbr\relax\else
\immediate\write16{Spaltenbreite hat sich geaendert!}\fi
\egroup}
\def\lcite#1{\edef\TTT{\noexpand\lowercase{\cite{#1}}}\TTT}
\Aussage{Example}
\def\tilde{\widetilde}
\def\rho{\varrho}

\def\hat{\widehat}
\def\Cinf{\cC^\infty}

\def\norm#1{|\!|#1|\!|}
\def\ri{{\it i)}}
\def\rii{{\it ii)}}
\def\riii{{\it iii)}}
\def\riv{{\it iv)}}

\fontdef Ss:cmss10,,.
\font\BF=cmbx10 scaled \magstep2
\font\CSC=cmcsc10 
\baselineskip15pt

{\baselineskip1.5\baselineskip\rightskip0pt plus 5truecm
\leavevmode\vskip0truecm\noindent
\BF Weyl Groups of Hamiltonian Manifolds, I

}
\vskip1truecm
\leftline{{\CSC Friedrich Knop}%
\footnote*{\rm Partially supported by a grant of the NSF}}
\leftline{Department of Mathematics, Rutgers University, New Brunswick NJ
08903, USA}
\leftline{knop@math.rutgers.edu}
\vskip1truecm
\beginsection Introduction. Introduction

Let $K$ be a compact connected Lie group and $M$ a compact Hamiltonian
$K$\_manifold, i.e., a symplectic $K$\_manifold equipped with a moment
map $\mu:M\pfeil\fk^*$. In this paper, we determine $\|Col|(M)$: the
set of all functions on $M$ which Poisson commute with all
$K$\_invariant functions. For this, we attach a finite reflection
group $W_M$ to $M$ and show that $\|Col|(M)$ is completely determined
by $\mu(M)$ and $W_M$.

More precisely, from these two data we construct a topological space
$Y$ equipped with a differentiable structure (in fact, it is
semi\_algebraic) and a (surjective) map $\hat\mu:M\pfeil Y$ such that
$\|Col|(M)$ consists exactly of the pull\_back functions via
$\hat\mu$. It is easy to see that, conversely, $\Cinf(M)^K$ is the
Poisson centralizer of $\|Col|(M)$. Thus we obtain a symplectic dual
pair
$$
\matrix{&&M\cr
&\swarrow&&\searrow\cr
Y&&&&M/K\cr}
$$
in the sense of Weinstein.

It follows immediately from the defining property of the moment map
$\mu$ that every pull\_back function via $\mu$ Poisson centralizes
$\Cinf(M)^K$. Thus, we obtain a homomorphism of Poisson algebras
$\mu^*:\Cinf(\fk^*)\pfeil\|Col|(M)$. This means that the moment map
$\mu$ factors through $Y$:
$$
\matrix{&&M\cr
&\Links{\hat\mu}\swarrow\cr
Y&&\Links\mu\downarrow\cr
&\Links\nu\searrow\cr
&&\fk^*\cr}
$$

In fact, Guillemin and Sternberg conjectured that $\mu^*$ is
surjective \cite{GS2}. This would mean that $\nu$ is a diffeomorphism
onto its image but a counterexample was given by Lerman in
\cite{Le0}. It turns out that, in his example, $Y$ is the (half-)cone
$x^2+y^2+z^2=t^2$, $t\ge0$ in $\RR^4$, $\fk^*=\RR^3$, and $\nu$ is the
projection to the first three coordinates. Then $t$ is not a
pull\_back of a {\it differentiable\/} function via $\nu$.  But $t$ is
the pull\_back of the {\it continuous\/} function
$\sqrt{x^2+y^2+z^2}$.

This is a general phenomenon: Karshon and Lerman, \cite{KL}, proved
that $\|Col|(M)=\mu^*\cC^0(\fk^*)\cap\Cinf(M)$. In our language, this
means that $\nu$ is a homeomorphism onto its image. This determines
$Y$ as a topological space and the complete determination of
$\|Col|(M)$ is a subtle problem of choosing a differentiable structure
on it.

This last step is controlled by the group $W_M$. It is a subquotient
of the Weyl group of $K$ and it is itself a reflection group. It
determines certain symmetry properties, the Taylor series of a
$\Cinf$\_function on $Y$ should have. Look, e.g., at Lerman's
counterexample above. There $W_M=\{1\}$. But in other situations also
$W_M=\{\pm1\}$ occurs. Then a function $f$ on $Y$ should be
differentiable if and only if $f$ is a differentiable function in
$x,y,z,t$ whose Taylor series at the vertex $(0,0,0,0)$ is
invariant under $(x,y,z,t)\mapsto(x,y,z,-t)$. This means that $f$
is a differentiable function of $(x,y,z)$ alone. The general case is
similar but more involved (see section \cite{Invariant collective} for
a complete statement.) 

The problem of determining $\|Col|(M)$ is local except that certain
connectivity properties of the fibers of $\mu$ are needed. Whenever
$M$ is compact, these hold by a theorem essentially due to Kirwan. But
we need to work locally. For that we axiomatize the properties needed
and introduce the notion of a {\it convex\/} Hamiltonian manifold. The
point is that every Hamiltonian manifold is locally convex. Then we
prove all results in this more general context. Note, that most
Hamiltonian manifolds appearing in practice (e.g., compact, or complex
algebraic, or cotangent bundle) are convex, as well.

Having localized the problem this way, we apply the symplectic slice
theorem. This allows us to assume that $M$ is an open subset of a {\it
real algebraic\/} Hamiltonian manifold $\Mq$. By powerful results of
Tougeron, or Bierstone-Milman, the property of being pull\_back can be
recognized on the level of Taylor series. This is now a purely
algebraic problem on $\Mq$. In particular, it can be solved by
complexifying $\Mq$. We show that the complexification of $\Mq$ is the
cotangent bundle $T^*_X$ of an affine $G$\_variety $X$ where $G$ is
the complexification of $K$ (thus reductive). In \cite{WuM}, I have
already determined all regular functions on $T^*_X$ which Poisson
commute with all $G$\_invariants.

Thus this paper consists mainly of two parts: sections \cite{Invariant
collective} and \cite{Proof of main theorem} deal with the comparison
of various classes of functions (differentiable$\leftrightarrow$power
series$\leftrightarrow$polynomials). In sections \cite{Alg Ham}
through \cite{Geometry II} we provide a crucial fact (\cite{Reduced})
about the geometry of $T_X^*$ which was not previously
established. Logically, these sections come first but since they are
of a very technical nature I postponed them.

For the convenience of the reader, two appendices are added. One
recalls comparison results between $\Cinf$\_functions and power
series. The other states the basic local structure theorems of
Hamiltonian varieties.

\medskip\noindent{\bf Acknowledgment:} I would like to thank Yael
Karshon and Eugene Lerman for explaining me their results in
\cite{KL}. I also benefited very much from discussions with Susan
Tolman and, in particular, Chris Woodward.

\medskip\noindent{\bf Notation:} Throughout this paper, ``analytic''
means ``real analytic''. Moreover, we use the term ``smooth'' always
in the algebraic geometric sense: ``smooth'' are morphism between
algebraic varieties which are surjective on tangent spaces. Smooth
functions in the sense of differential geometry are called $\Cinf$ or
differentiable. Moreover, differentiable always means infinitely often
differentiable.

\beginsection Topology. Convex Hamiltonian manifolds

Let $K$ be a connected, compact Lie group with Lie algebra $\fk^*$ and
$M$ a Hamiltonian $K$\_manifold with moment map $\mu:M\pfeil\fk^*$.

Let $\ft\subseteq\fk$ be a Cartan subalgebra corresponding to a
maximal torus $T\subseteq K$ and $W=W(\fk,\ft)$ its Weyl group. Since
$\ft$ has a unique $T$\_stable complement in $\fk$, we can regard
$\ft^*$ as a subspace of $\fk^*$. The map $\ft^*/W\pfeil\fk^*/K$
is a diffeomorphism. Furthermore, its restriction to a Weyl chamber
$\ft^*_+\subseteq\ft^*$ is a differentiable homeomorphism, but in
general not a diffeomorphism. We use it to construct a
continuous map $\psi$ which makes the following diagram commutative:
$$
\matrix{
M&\pf\mu&\fk^*\cr
\Links\psi\untenPf&&\untenPf\cr
\ft^*_+&\pfeil&\ft^*/W\cr}
$$
The map $\psi$ is also characterized by the property
$K\mu(x)\cap\ft^*_+=\{\psi(x)\}$ for all $x\in M$.

In this section we study some topological properties of $\psi$. For any two
(not necessarily distinct) points $u,v\in\ft^*$ let $[u,v]$ be the
line segment joining them.

\Definition: A Hamiltonian $K$\_manifold is called {\it convex\/} if
$\psi^{-1}([u,v])$ is connected for all $u,v\in\psi(M)$.
\medskip\noindent
Clearly, $M$ is convex if and only if $\psi^{-1}(B)$ is connected for
every convex subset $B$ of $\ft^*$. In practice we will use another
characterization of convexity:

\Theorem ConvexCrit. A Hamiltonian $K$\_manifold is convex if and only
if the following conditions are satisfied:
\item\ri The fibers of $\psi$ are connected.
\item\rii The image $\psi(M)$ is convex.
\item\riii The map $\psi:M\pfeil\psi(M)$ is open.\Par

\Proof: Assume first that \ri{} through \riii{} hold. If
$\psi^{-1}([u,v])$ is disconnected there are $U_1,U_2\subseteq M$ open
such that $\psi^{-1}([u,v])$ is the disjoint union of $X_1$ and $X_2$
where $X_i:=\psi^{-1}([u,v])\cap U_i\ne\emptyset$. By \ri{} we have
$\psi(X_1)\cap\psi(X_2)=\emptyset$ and, by \rii,
$[u,v]=\psi(X_1)\cup\psi(X_2)$. Finally, \riii{} implies
that $\psi(X_i)=[u,v]\cap\psi(U_i)$ is open in $[u,v]$. This
contradicts the connectedness of $[u,v]$.

For the reverse direction, we need some preparation. Let $V$ be a real
vector space equipped with a lattice $\Gamma\subseteq V$. Our example
will be $V=\ft^*$ with the lattice consisting of all ${1\over i}d\chi$
where $T\pfeil\CC^*$ is a character. A {\it homogeneous rational
cone\/} is a subset of the form $\sum_{i=1}^N\RR^{\ge0}\gamma_i$ where
$\gamma_1,\ldots,\gamma_N\in\Gamma$. A {\it rational cone\/} is a
translate $u+C$ of a homogeneous rational cone $C$ by a vector $u\in
V$. In this case we say that $u$ is a vertex of $u+C$. The following
theorem of Sjamaar is the foundation for the results in this section.

\Theorem SjaThm. {\rm(\cite{Sja}~Thm.~5.5)} Let $M$ be a Hamiltonian
$K$\_manifold. Then for every orbit $Kx\subseteq M$ there is a unique
rational cone $C_x\subseteq\ft^*$ with vertex $\psi(x)$ such that:
\item{\ri} There exist arbitrary small $K$\_stable neighborhoods $U$
of $Kx$ such that $\psi(U)$ is a neighborhood of $\psi(x)$ in $C_x$.
\item{\rii} For $u\in\ft^*_+$ let $x$ and $y$ be in the same connected
component of $\psi^{-1}(u)$. Then $C_x=C_y$.\Par

Sjamaar's theorem implies a general openness property of $\psi$:

\Lemma Open. Let $M$ be a Hamiltonian $K$\_manifold such that all
fibers of $\psi:M\pfeil\ft^*_+$ are connected. Let $U\subseteq M$ be
open. Then also $\psi^{-1}\psi(U)$ is open.

\Proof: We may assume that $U$ is $K$\_stable. For every
$y\in\psi^{-1}\psi(U)$ there is $x\in U$ such that
$\psi(x)=\psi(y)=:u$. Since $\psi^{-1}(u)$ is connected,
\cite{SjaThm}{\rii} implies $C_x=C_y=:C$. By part \ri{} of that
theorem there are open neighborhoods $U_x$, $U_y$ of $x$, $y$,
respectively, such that $\psi(U_x)$ and $\psi(U_y)$ are neighborhoods
of $u$ in $C$. Hence $(\psi|_{U_y})^{-1}(\psi(U_x))$ is a neighborhood
of $y$ which is contained in $\psi^{-1}\psi(U)$.\qed

\Lemma Conn. Let $M$ be a Hamiltonian $K$\_manifold such that all
fibers of $\psi:M\pfeil\ft^*_+$ are connected. Assume moreover that
every $u\in\psi(M)$ has arbitrary small neighborhoods $B$ such that
$\psi^{-1}(B)$ is connected. Then $\psi:M\pfeil\psi(M)$ is an open
map.

\Proof: Let $x\in M$ and $U$ an open neighborhood of $x$. We have to
show that $\psi(U)$ is a neighborhood of $u:=\psi(x)$ in
$\psi(M)$. For this we may assume that $\psi(U)$ is a neighborhood of
$u$ in $C_x$. Let $B$ be a neighborhood of $u$ in $\ft^*_+$ such that
$B\cap C_x\subseteq\psi(U)$ and such that $\psi^{-1}(B)$ is
connected. Let $V_1:=\psi^{-1}\psi(U)$ which is open in $M$ by
\cite{Open}. Clearly, also $V_2:=\psi^{-1}(\ft^*\setminus C_x)$ is
open in $M$. Moreover, $V_1$ and $V_2$ are disjoint and cover
$\psi^{-1}(B)$. Connectivity implies $\psi^{-1}(B)\subseteq V_1$,
i.e., $B\cap\psi(M)\subseteq\psi(U)$ which proves the assertion.\qed

Now we can complete the proof of \cite{ConvexCrit}. Assume that $M$ is
convex. Then \ri{} and \rii{} clearly hold. Let $u\in\psi(M)$ and
$B\subseteq\ft^*_+$ a convex neighborhood of $x$. Since also $\psi(M)$
is convex we have $[u,v]\subseteq B\cap\psi(M)$ for every $v\in
B\cap\psi(M)$. By assumption, $\psi^{-1}([u,v])$ is connected. This
implies that $\psi^{-1}(B)$ is connected. Thus \riii{} holds by
\cite{Conn}.\qed

\Remark: I don't know of any example of a Hamiltonian $K$\_manifold $M$
where $\psi(M)$ is convex, the fibers of $\psi$ are connected but
where $\psi$ is not an open map.
\medskip
If $M$ is convex then $C_x$ depends by \cite{SjaThm} only on
$u=\psi(x)$. Thus we may set $C_u:=C_x$. Let $C$ be a rational cone
with vertex $u$. For $v\in C$ let $T_vC:=C+\RR^{\ge0}(u-v)$ be the
{\it tangent cone\/} of $C$ in $v$. This is the smallest rational cone
containing $C$ and having $v$ as a vertex.

\Theorem ConvProp. Let $M$ be a convex Hamiltonian $K$\_manifold and
$u\in\psi(M)$.
\item{\ri} The image $\psi(M)$ is contained in $C_u$. Moreover, it is
a neighborhood of $u$ in $C_u$.
\item{\rii}The image $\psi(M)$ is locally a polyhedral cone and, in
particular, locally closed and semi\_analytic.
\item{\riii} For every $v\in\psi(M)$ in a neighborhood of $u$ holds
$C_v=T_vC_u$.

\Proof: \ri{} Let $x\in M$ with $u=\psi(x)$. Then $Kx$ has an open
neighborhood $U$ such that $\psi(U)$ is a neighborhood of $u$ in
$C_u$. Let $v\in\psi(M)$. Since $\psi(U)$ is open in $\psi(M)$ and
since $[u,v]\subseteq\psi(M)$ also $[u,v]\cap\psi(U)$ is open in
$[u,v]$. This implies $[u,v]\subseteq C_u$, thus $\psi(M)\subseteq
C_u$. Moreover, $\psi(M)$ is a neighborhood of $u$ in $C_u$ since
already $\psi(U)$ is.

\rii{} Follows directly from \ri.

\riii{} By \ri{} there is an open neighborhood $U$ of $x$ such that
$\psi(U)$ is open in $C_u$. Let $v\in\psi(U)$. Then \ri, applied to
$v$, implies that $C_v$ is the cone spanned by $\psi(U)$ over the
vertex $v$. On the other hand, $\psi(U)$ is open in $C_u$. Thus
$T_vC_u$ is also the cone spanned by $\psi(U)$ over $v$.\qed

Let $V$ be a unitary representation of $K$. Then every smooth
$K$\_stable complex algebraic subvariety of ${\bf P}(V)$ is in a
canonical way a Hamiltonian $K$\_manifold. We call a Hamiltonian
$K$\_manifold {\it projective} if it arises this way but possibly with
the symplectic form and the moment map rescaled by some non\_zero
factor.

\Proposition Proj. Let $M$ be a Hamiltonian $K$\_manifold such that for
every projective Hamiltonian $K$\_manifold $X$ holds that $\psi_{M\times
X}^{-1}(0)$ is connected. Then $\psi_M:M\pfeil\psi(M)$ is an open map
with connected fibers.

\Proof: Let $\Xq$ equal $X$ with the symplectic structure multiplied
by $-1$. Then we have $\mu_{M\times\Xq}(m,x)=\mu_M(m)-\mu_X(x)$. Thus
choosing for $X$ the coadjoint orbit $Ku$, $u\in\ft^*_+$ we obtain
that $\psi_M^{-1}(u)=\psi_{M\times\Xq}^{-1}(0)$ is connected. Now let
$X_0$ be projective such that $\psi_{X_0}(X_0)$ is a neighborhood of $0$,
e.g., $X_0={\bf P}(V)$ where $V$ contains stable points for the
$K^\CC$\_action. By rescaling, we can arrange that $\psi(X_0)$ is
arbitrary small. Let $X:=X_0\times Ku$. Then $B:=\psi_X(X)$ is an
arbitrary small neighborhood of $u$ in $\ft^*_+$. Moreover, projection
to the first factor induces a surjective map
$\psi_{M\times\Xq}^{-1}(0)\auf\psi_M^{-1}(B)$. By assumption the
first, hence the second set is connected. We conclude with
\cite{Conn}, that $\psi_M$ is open onto its image.\qed

\noindent Now we can give examples of convex Hamiltonian manifolds:

\Theorem Examp. Every Hamiltonian $K$\_manifold $M$
satisfying one of the conditions {\ri}, {\rii}, {\riii} below is
convex.
\item{\ri} The moment map $\mu:M\pfeil\fk^*$ is proper (e.g., if $M$ is
compact).
\item{\rii} The manifold $M$ is a complex algebraic variety, the action
of $K$ is the restriction of an algebraic $K^\CC$\_action, and the symplectic
structure is induced by a $K$\_invariant K{\"a}hler metric.
\item{\riii} The manifold $M$ is a complex Stein space, the action
of $K$ is the restriction of a holomorphic $K^\CC$\_action, and the symplectic
structure is induced by a $K$\_invariant K{\"a}hler metric.\Par

\Proof: In all cases, it is known that $\psi$ has connected fibers and
convex image (see \cite{HNP} or \cite{Sja} for {\ri} and
\cite{HeHu} for {\rii} and {\riii}). Moreover, the classes {\ri}-{\riii} are
preserved by taking the product with a projective Hamiltonian
$K$\_manifold. Thus \cite{Proj} implies that $\psi$ is open.\qed

\noindent Locally, every Hamiltonian manifold is convex:

\Theorem SKL. Let $M$ be a Hamiltonian $K$\_manifold. Then every $x\in
M$ has a convex $K$\_stable open neighborhood $U$ such that $\psi(U)$
is open in $C_x$.

\Proof: The proof is similar to that of \cite{KL} Prop.~3.7. First, if
the orbit $Kx$ (with $\mu(x)\in\ft^*_+$) is not isotropic then let
$L=K_{\mu(x)}$ and replace $M$ by $M_L\cong K\times^L M(L)$
(\cite{Cross}). If there is already a neighborhood $U(L)$ of $x$ in
$M(L)$ with the properties claimed in the theorem with respect to $L$,
then it is trivial to check that $U:=K\cdot U(L)$ has the required
properties for $K$.

Thus we may assume that $Kx$ is isotropic. By \cite{Slice} and the
remark following it, we may assume that $M=K\times^H(\fh^\perp\times
S)$ and $x=[1,(0,0)]$. Here, $H$ is a closed subgroup of $K$ and $S$
is a unitary representation of $H$.

Let $D$ and $\Dq$ be the open and the closed unit ball in $S$. Put
$V:=K\times^H(\fh^\perp\times D)$,
$\Vq:=K\times^H(\fh^\perp\times\Dq)$. The restriction of $\mu$ to
$\Vq$ factors as
$$
\mu|_\Vq:\Vq\into K\times^H(\fk^*\times\Dq)\cong
(K\times^H\Dq)\times\fk^*\auf\fk^*.
$$
This shows that $\mu|_\Vq$ is proper.

We let $t\in\RR^{>0}$ act on $M$ by
$t\cdot[k,(u,v)]:=[k,(t^2u,tv)]$. Then $\psi$ becomes homogeneous,
more precisely, $\psi(t\cdot y)=t^2\psi(y)$ for all $y\in M$.

There is a $K$\_stable open neighborhood $U_0$ of $x$ in $V$ such that
$\psi(U_0)$ is a neighborhood of $0=\psi(x)$ in $C_x$. From
$M=\RR^{>0}\cdot U_0$ we obtain $\psi(M)\subseteq C_x$. In particular,
also $\psi(V)$ is a neighborhood of $0$ in $C_x$. Thus there is a
convex, open neighborhood $B$ of $0$ in $\ft^*$ such that $B\cap
C_x\subseteq\psi(V)$ and $\Bq\cap C_x\subseteq\psi(\Vq)$. We put
$U:=V\cap\psi^{-1}(B)$.

By construction, $\psi(U)=B\cap C_x$ is convex and open in $C_x$.  It
remains to show that $\psi: U\pfeil\psi(U)$ is open with connected
fibers. The closure of $U$ is $\Uq=\Vq\cap\psi^{-1}(\Bq)$. Moreover,
we have $U=\bigcup_{0<t<1}t\cdot\Uq$. Since increasing unions of
connected sets are connected and by homogeneity it suffices to show
that $\psi:\Uq\pfeil\Bq$ has connected fibers. We show that even
$\psi:\Vq\pfeil C_x$ has connected fibers.

For this recall that $S$ is a unitary vector space. In particular,
$U(1)$ acts on $S$ and thus on $M$. This action is Hamiltonian with
moment map $\mu_0([k,(u,v)])=\norm v^2$. Let $P$ be the symplectic
cut of $M$ at the level $\mu_0=1-0$: As a set we have
$P:=P_1{\buildrel.\over\cup}P_2$ where $P_1:=\mu_0^{-1}([0,1))$ and
$P_2:=\mu_0^{-1}(1)/U(1)$. It turns out that $P$ is a symplectic
manifold, see \cite{Le} for details.

There is an obvious proper map $\pi:\Vq\auf P$ where $V$ goes
isomorphically onto $P_1$ and the $U(1)$\_orbits in the boundary are
collapsed to points in $P_2$. Since the $U(1)$\_action commutes with
the $K$\_action, $P$ is even a Hamiltonian $K$\_manifold. Its moment
map is the unique map which makes the following diagram commutative:
$$
\matrix{
\Vq\cr
\Links\pi\downarrow&\searrow\Rechts\mu\cr
P&\pf{\mu_P}&\fk^*\cr}
$$
In particular, $\mu_P$ is proper. This implies that $\mu_P$, and
$\psi_P$ have connected fibers (\cite{HNP}, \cite{Sja}~1.4). Since
$\pi$ is proper with connected fibers, also the fibers of $\psi|_\Vq$
are connected, which was claimed.

By \cite{Examp}\ri, the map $\psi_P:P\pfeil\psi(P)$ is open. Since
$U\subseteq P$ is open, this implies that $\psi:U\pfeil\psi(U)$ is
open, as well.\qed

An application is the following well known statement.

\Corollary fa0. Let $M$ be any Hamiltonian $K$\_manifold and let
$\fa^0\subseteq\ft^*$ be the affine subspace spanned by
$\psi(M)$. Then $\fa^0$ is also the affine span of every local cone
$C_x$, $x\in M$. Moreover, $\psi(M)$ is Zariski\_dense in $\fa^0$.

\Proof: With local convexity and \cite{ConvProp}\riii{} it follows that
the affine span of $C_x$ is locally constant, hence constant in
$x$.\qed

The last class of examples for convex Hamiltonian manifold was
suggested to me by S.~Tolman:

\Theorem. Let $X$ be a $K$\_manifold and $M:=T^*_X$ with its natural
Hamiltonian structure. Then $M$ is convex.

\Proof: The map $\psi:M\pfeil\ft^*$ is homogeneous for the natural
scalar action on the fibers. Furthermore, $X$ is embedded into $M$ as
the zero\_section with $X\subseteq\psi^{-1}(0)$. Therefore, the cone
$C=C_x$ does not depend on $x\in X$. Thus, $X$ has a neighborhood $U$
such that $\psi(U)$ is a neighborhood of $0$ in $C$. Since $\psi$ is
homogeneous, we conclude $\psi(M)=C$. This shows in particular that
$\psi(M)$ is convex.

Let $\pi:M\pfeil X$ be the projection. Let $m\in M$ and
$x=\pi(m)$. Then, again by homogeneity, there is a convex neighborhood
$V$ of $Kx$ such that $m\in V$ and such that $\psi(V)$ is open in $C$
(\cite{SKL}). Then any neighborhood of $m$ is mapped to a neighborhood
of $\psi(m)$ in $\psi(V)$. This implies that $\psi:M\pfeil\psi(M)$ is
open.

Finally, assume that the fiber over $u\in\psi(M)$ is not connected,
i.e., is a disjoint union of two non\_empty open pieces $F_1$ and
$F_2$. For every open $K$\_stable subset $U$ of $X$ holds
$\psi(\pi^{-1}(U))=C$. Hence, $\pi(\psi^{-1}(u))$ is dense in
$X$. Thus there is
$x\in\overline{\pi(F_1)}\cap\overline{\pi(F_2)}$. Let $V$ be a convex
neighborhood of $Kx$ in $M$. By homogeneity, we may assume $V$ meets
both $F_1$ and $F_2$. But this contradicts the connectedness of fibers
of $\psi|_V$. Thus, $\psi$ has connected fibers.\qed

\beginsection Invariant collective. Collective functions: Statement of
the result

Let $M$ be a Hamiltonian $K$\_variety. Every element $\xi\in\fk$
induces a function $\mu_\xi$ in $M$ by $\mu_\xi(x):=\<\mu(x),\xi\>$.
One of the defining properties of a moment map is
$$
\{f,\mu_\xi\}=df(\xi_*)\quad\hbox{for all $f\in\Cinf(M)$.}
$$
Here $\xi_*$ is the vector field expressing the infinitesimal action
of $\xi$ on $M$. Thus, since $K$ is connected, a function $f$ is
$K$\_invariant if and only if it Poisson commutes with all functions
$\mu_\xi$. Conversely, consider the set of {\it collective functions\/}
$$
\|Col|(M):=\big\{g\in\Cinf(M)\mid\{\Cinf(M)^K,g\}=0\,\big\}.
$$
Thus every $\mu_\xi$ is collective. Slightly more generally,
pull\_back by $\mu$ induces a homomorphism of Poisson algebras
$\mu^*:\Cinf(\fk^*)\pfeil\|Col|(M)$.

The Poisson subalgebras $\Cinf(M)^K$ and $\|Col|(M)$ form a ``dual
pair'' in $\Cinf(M)$, i.e., one is the Poisson centralizer of the
other. Their intersection, the set of invariant collective functions
$\|Col|^K(M)$, is the Poisson center of either algebra and is
particularly important. Note that $\mu^*$ induces a homomorphism
$\Cinf(\ft^*/W)=\Cinf(\fk^*)^K\pfeil\|Col|^K(M)$.

Let $U\subseteq M$ be open, $K$\_stable. Then
$\|res|^M_U:\Cinf(M)\pfeil\Cinf(U)$ has a dense image which implies that
the restriction of a collective function to $U$ is again collective.

Our aim is to describe the set of all collective functions. One step
in this direction is the following criterion due to Karshon \cite{Ka}
(see \cite{KL}~Thm.~1).

\Theorem KarLer. Let $M$ be a Hamiltonian $K$\_manifold. Then $f\in\Cinf(M)$
is collective if and only if $f$ is locally constant on the fibers of
$\mu$.

This theorem shows in particular that being collective is mostly a
set\_theoretic property. In particular, if a sequence of collective
functions converges pointwise to a differentiable function $f$ then also $f$
is collective.

Our goal is a more precise description of $\|Col|(M)$. For this we
introduce first a basic construction. Let again $\psi:M\pfeil\ft^*_+$
be the map with $K\mu(x)\cap\ft^*_+=\{\psi(x)\}$. Let
$\fa^0\subseteq\ft^*$ be the affine subspace spanned by $\psi(M)$ (see
\cite{fa0}). We let $K$ act on $V:=\fk^*\times\fa^0$ by conjugation
on the first and trivially on the second factor and consider
$K\fa^0\subseteq V$ where $\fa^0$ is embedded diagonally into
$V$. This is a closed subset. In fact it equals (as a set) the fiber
product $\fk^*\times_{\ft^*/W}\fa^0$, i.e., the set of pairs $(u,v)$
such that $u$ and $v$ have the same image in $\fk^*/K=\ft^*/W$. This
shows moreover that $K\fa^0$ is (the set of real points of) a
real\_algebraic closed subvariety of $V$. We define $\RR[K\fa^0]$ to
be the image of $\RR[V]\pf{\|res|}\Cinf(K\fa^0)$. Note, that the
kernel of $\RR[V]\auf\RR[K\fa^0]$ is in general (much) larger than the
ideal generated by the obvious equations expressing $K\fa^0$ as a
fiber product. For example, the most extreme case would be $\fa^0=0$,
where the fiber product is the nilcone of $\fk^*$ which happens to
have just one real point but is nevertheless of positive dimension as
an algebraic variety.

\Definition: Let $\tilde{K\fa^0}$ be (the set of real points of) the
normalization of $K\fa^0$, i.e., $\tilde{K\fa^0}=\|AlgHom|(A,\RR)$
where $A$ is the integral closure of $\RR[K\fa^0]$ in its field of
fractions. Moreover, we equip $\tilde{K\fa^0}$ with the real
(Hausdorff) topology.
\medskip

\noindent By definition, there is a natural morphism
$\nu:\tilde{K\fa^0}\pfeil K\fa^0$ and the $K$\_action lifts to
$\tilde{K\fa^0}$.

\Proposition Quotient. The map $\nu:\tilde{K\fa^0}\pfeil K\fa^0$ is a
homeomorphism and $\tilde{K\fa^0}/K\pfeil\fa^0$ is a diffeomorphism.

\Proof: Observe
$(\fk^*\times_{\ft^*/W}\fa^0)/K=\fk^*/K\times_{\ft^*/W}\fa^0=\fa^0$.
Thus $K\fa^0/K=\fa^0$. Since $\fa^0$ is already normal we have also
$\tilde{K\fa^0}/K=\fa^0$. In particular, $\nu$ induces an isomorphism
on orbit spaces. The fibers of $\nu$ are finite and all isotropy
groups in $K\fa^0$ (being Levi subgroups) are connected which implies
that $\nu$ is injective. Moreover, the dense subset of regular points
of $K\fa^0$ is contained in the image of $\nu$. Since $\nu$ is proper,
it must be a homeomorphism.\qed

\noindent Thus $\tilde{K\fa^0}$ equals $K\fa^0$ as a topological space
but it may have more differentiable functions. The image of
$\mu\times\psi:M\pfeil V$ is clearly contained in $K\fa^0$. So we
obtain a map $\tilde\mu:M\pfeil\tilde{K\fa^0}$ which is continuous but
in general not differentiable.

For studying the local structure of $\tilde{K\fa^0}$ the following
lemma is helpful. For a Levi subgroup $L$ let $\ft^r\subseteq\ft^*$ be
the set of points $v$ with $W_v\subseteq W_L$. Put $\fa^r:=\fa^0\cap
\ft^r$ (which might be empty).

\Lemma Scheibe1. For any Levi subgroup $L$ there is a canonical map
$\lambda:K\times^L\tilde{L\fa^0}\pfeil\tilde{K\fa^0}$
which is a diffeomorphism over $\fa^r$.

\Proof: There is clearly a surjective map $K\times^LL\fa^0\pfeil
K\fa^0$. Then $\lambda$ is just the unique lift to the normalizations.
To show that it is an isomorphism over $\fa^r$ we complexify. Let
$G:=K_\CC$. Then $(K\fa^0)_\CC$ is the Zariski closure
$\overline{G\fa^0_\CC}$. \cite{Scheibe} implies that
$G\times^{L_\CC}\overline{L_\CC\fa^0_\CC}\pfeil\overline{G\fa^0_\CC}$
is a closed embedding over $\fa^r_\CC$. Since it is clearly surjective
it is an isomorphism over $\fa^r_\CC$. Hence the same holds for the
normalizations.\qed

The difference between $\tilde{K\fa^0}$ and $K\fa^0$ is not yet
completely understood. It translates into a property of certain {\it
nilpotent\/} orbits of the {\it complexified\/} Lie algebra. We will
state two results in this direction which might be useful. Let's call
a compact group unitary if it is locally isomorphic to a product of a
torus and a product of special unitary groups.

\Proposition SU(n)-1. \ri{} Assume $\fa^0$ contains an interior point of
the Weyl chamber $\ft^*_+$. Then $\nu$ is a diffeomorphism.\Par\noindent
\rii{} Let $v\in\fa^0$ and assume that the
isotropy group $K_v$ is unitary (this holds, in particular, if $K$
itself is unitary). Then $\nu$ is a diffeomorphism in a
neighborhood of $v$.

\Proof: In case \ri{} it suffices to show that $\nu$ is a
diffeomorphism near every $v\in\fa^0$. Then with \cite{Scheibe1} we
may reduce to the case $v=0$. Let $L\subseteq K$ be the centralizer of
$\fa^0$. Then either $L=T$ or $K$ is unitary. We are going to show
that then $\nu$ is even an isomorphism of algebraic varieties.

For this we may complexify. Then $L_\CC\subseteq G:=K_\CC$ is the Levi
subgroup of a parabolic $P$. Let $P_u$ be its unipotent radical. The
complexification of $K\fa^0$ is $X:=\overline{G\fa^0_\CC}$. By the
choice of $L$ we have $\overline{P\fa^0}=\fr:=\fa^0\oplus\fp_u$. Thus
the proper morphism $G\times^P\fr\pfeil\fg$ has exactly image $X$. It
factors through the normalization $\XS$ of $X$. Thus we have
$G\times^P\fr\pf\pi\XS\pf\nu X$. These are morphisms over $\fa^0$.

Since $L$ is the centralizer of $\fa^0$, the morphism $\pi$ is
birational. Then it follows from \cite{WuM}~4.1 that $\XS$ has
rational singularities and is in particular Cohen\_Macaulay. Let
$\XS_0$ be the zero\_fiber of $\XS\pfeil\fa^0$. Since $\pi$ is an
isomorphism over a generic point of $\XS_0$ it is generically
reduced. Being defined by a regular sequence it is reduced. By
Nakayama's lemma, $\nu$ is an isomorphism if and only if its schematic
zero\_fiber $\nu^{-1}(0)$ is a (reduced) point. Since
$\nu^{-1}(0)\subseteq\XS_0$ we are reduced to show that
$\XS_0\pfeil\fg$ is a closed embedding. This would we true if $G\fp_u$
is normal and $G\times^P\fp_u\pfeil G\fp_u$ is birational. But these
conditions hold under our assumptions: if $L=T$ then $P=B$ is a Borel
subgroup and we conclude by results of Kostant and Springer (Springer
resolution, see e.g.~\cite{Spr}). If $K$ is unitary then all nilpotent
orbit closures are normal (Kraft\_Procesi~\cite{KP}) and all
stabilizers are connected modulo center.\qed

The affine subspace $\fa^0$ of $\ft^*$ is in general not
$W$\_stable. Thus we define
$$
\eqalign{
N_W(\fa^0)&:=\{w\in W\mid w\fa^0=\fa^0\};\cr
C_W(\fa^0)&:=\{w\in W\mid w|_{\fa^0}=\|id|\};\cr
W(\fa^0)&:=N_W(\fa^0)/C_W(\fa^0).\cr}
$$
If we let $W(\fa^0)$ act on the second factor of $V=\fk^*\times\fa^0$
we obtain actions on $K\fa^0$ and $\tilde{K\fa^0}$ which commute with
$K$. Let $W'\subseteq W(\fa^0)$ be any subgroup. Then also
$\tilde{K\fa^0}/W'\pfeil K\fa^0/W'$ is a homeomorphism and
$(\tilde{K\fa^0}/W')/K=\fa^0/W'$. Now we can formulate our main
result.

\Theorem Main. Let $M$ be a convex Hamiltonian $K$\_manifold and let
$\fa^0$ be the affine subspace spanned by $\psi(M)$. Then there is
a unique subgroup $W_M$ of $W(\fa^0)$ such that:
\item{\ri} The map $\Phi=\tilde\mu/W_M:M\pfeil\tilde{K\fa^0}/W_M$ is
differentiable and induces an isomorphism $\Cinf(\Phi(M))\pf\sim\|Col|(M)$.
\item{\rii} Let $R$ be the set of all reflections in $W(\fa^0)$ whose
reflecting hyperplane meets $\psi(M)$. Then $W_M$ is generated by a
subset of $R$.\Par

\noindent The proof is given in the next section. Let me first give
some conclusions.

\Corollary. The map $\phi=\psi/W_M:M\pfeil\fa^0/W_M$ is differentiable and
induces an isomorphism $\Cinf(\phi(M))\pf\sim\|Col|^K(M)$. Thus
$\|Col|^K(M)$ can be identified with the algebra of differentiable functions
on an $r$\_dimensional semi\_analytic subset of $\RR^r$ where
$r=\|dim|\fa^0$.

\Proof: The first statement follows from \cite{Quotient} and
\cite{Main}\ri{} by taking $K$\_invariants. The second follows from
the fact that the ring of invariants of a reflection group is a
polynomial ring (hence $\fa^0/W_M\into\RR^r$) and that $\psi(M)$ is
locally even semi\_algebraic (\cite{ConvProp}\rii).\qed

\Corollary SU(n)-2. \ri{} The ring $\|Col|(M)$ is a finitely generated
$\Cinf(\mu(M))$\_module.\Par\noindent
\rii{} Assume that $\psi(M)$ contains an interior point of $\ft^*_+$ or
that $K_v$ is unitary for all $v\in\psi(M)$. Then
$\|Col|(M)$ is, as a $\Cinf(\mu(M))$\_module, generated by
finitely many {\rm invariant\/} collective functions.

\Proof: \ri{} The morphism of algebraic varieties
$\tilde{K\fa^0}\pfeil\fk^*$ is finite. Hence there are
$f_1,\ldots,f_s$ in $\RR[\tilde{K\fa^0}]^{W_M}$ which generate it as a
$\RR[\fk^*]$\_module. It follows from \cite{Mal} Ch.~V, Cor.~4.4 that
for every $x\in M$ the $f_i$ generate the stalk
$\Cinf_{\Phi(M),\Phi(x)}$ as a $\Cinf_{\mu(M),\mu(x)}$\_module. By a
partition of unity argument this shows that $\|Col|(M)=\Cinf(\Phi(M))$
is generated by the $f_i$ as a $\Cinf(\mu(M))$\_module.

\rii{} The assumptions imply (\cite{SU(n)-1}) that
$\tilde{K\fa^0}\pfeil K\fa^0\subseteq\fk^*\times\fa^0$ is a
diffeomorphism in a neighborhood of $\tilde\mu(M)$. Thus we can
replace $\tilde{K\fa^0}$ by $K\fa^0$ and therefore choose the $f_i$ to
be $K$\_invariant.\qed

\noindent Another characterization of $W_M$ is:

\Proposition Minimal. A subgroup $W'$ of $W(\fa^0)$ contains $W_M$ if
and only if $\psi/W':M\pfeil \fa^0/W'$ is differentiable.

\Proof: If $W_M\subseteq W'$ then $\phi':=\psi/W'$ is clearly
differentiable. Conversely, assume $\phi'$ is differentiable. Then \cite{KarLer}
implies $\Cinf(\phi'(M))\into\|Col|^K(M)=\Cinf(\phi(M))$. This shows
that $\phi(M)\pfeil\phi'(M)$ is differentiable. Looking at the formal stalk at
$u\in\psi(M)$ this implies $\Cinf(\psi(M))^{\wedge{W'_u}}_u\subseteq
\Cinf(\psi(M))^{\wedge{W_{M,u}}}_u$ and therefore $W_{M,u}\subseteq
W'_u$. Hence $W_M\subseteq W'$ since $W_M$ is generated its
isotropy groups in $\psi(M)$.\qed

We conclude this section with two examples. The first one is just
Lerman's example from \cite{Le0}: let $K=SU(2)$ and $M:=\CC^2$ with
its symplectic structure coming from the $K$\_invariant Hermitian form
$q(z_1,z_2):=z_1\zq_1+z_2\zq_2$. Then a moment map $M\pfeil\fk^*$
exists and is in fact given by {\it quadratic\/} polynomials. Every
$K$\_invariant on $M$ is a composite with $q$. Hence, $q$ is invariant
collective. Assume there is $h\in\Cinf(\fk^*)$ with
$q=f\circ\mu$. Then $h$ were a non\_zero $K$\_invariant function of
degree one which clearly does not exist. On the other hand, one
easily verifies that $\psi:M\pfeil\ft^*\cong\RR$ is given by
$q$. Thus $\psi$ is differentiable and $W_M=1$.

The second example is a variant of the preceding one: Let
$M:=\CC^{2n}$ and $K$ the maximal compact subgroup of
$Sp_{2n}(\CC)$. Then again a quadratic moment map $\mu:M\pfeil\fk^*$
exists and $K$ fixes a Hermitian norm $q$. Let $\Kq:=U(2n)$ be the
unitary group associated to $q$. Then there is also a quadratic moment
map $\overline\mu:M\pfeil\overline\fk^*$ such that $\mu$ is the
composition of $\overline\mu$ with $\overline\fk^*\pfeil\fk^*$.  It
happens that every $K$\_invariant is a composite with $q$. Thus
$\Cinf(M)^K=\Cinf(M)^\Kq$ and therefore,
$\|Col|_K(M)=\|Col|_\Kq(M)$. This shows in particular that
$\overline\fk^*$ embeds into $\|Col|_K(M)$ as a space of quadratic
polynomials. But the space of quadratic polynomials which pull\_back
from $\mu\times\psi$ has at most dimension $\|dim|K+1$. This shows
that for $n>1$ the normalization of $K\fa^0$ is really necessary. It
also shows that one needs some restrictions in \cite{SU(n)-2}.

\beginsection Proof of main theorem. Proof of \cite{Main}

We investigate the problem of finding the collective functions first
on the level of formal power series. For this, we adopt the following
notation: If $X$ is a set, $Y\subseteq X$ a subset and $\cR$ a ring of
functions on $X$, then let $\cR^\wedge_Y$ denote the completion of
$\cR$ with respect to the ideal of functions vanishing in $Y$.

Assume that $Kx$ is isotropic orbit of $M$. By \cite{Slice} and the
remark following it we may assume that $M=K\times^H(\fh^\perp\times
S)$ and $x=[1,(0,0)]$. The manifold $M$ has a real algebraic
structure: with $\RR[K]\subseteq\Cinf(K)$, the ring of representative
(i.e., $K$\_finite) functions, and $\RR[\fh^\perp\times S]$, the ring
of real polynomials, we define
$$
\RR[M]:=(\RR[K]\otimes\RR[\fh^\perp\times S])^H.
$$
This is a finitely generated $\RR$\_algebra and turns $M$ into a real
algebraic manifold. All $K$\_orbits are Zariski closed. Also the symplectic
structure and the moment map are real algebraic. Thus $\RR[M]$ is a Poisson
algebra, $\RR[M]^K$ a Poisson subalgebra. Let $\|Col|_{\rm
alg}(M)\subseteq\RR[M]$ be its centralizer.

\Lemma Comparison. Let $M$ and $x$ be as above. Then there is a
subgroup $W(x)$ of $W(\fa^0)$ which is generated by reflections such
that $\Phi=\tilde\mu/W(x):M\pfeil\tilde{K\fa^0}/W(x)$ is algebraic and
induces an isomorphism
$\RR[\tilde{K\fa^0}/W(x)]^\wedge_0\pf\sim\|Col|_{\rm alg}(M)^\wedge_{Kx}$.

\Proof: We are going to determine the complexification of the
real\_algebraic variety $M$. Since $H$ is compact, the symplectic
structure of $S$ comes from an $H$\_stable Hermitian form. In
particular, $S$ has also the structure of a complex vector
space. Hence $S_\CC=S\otimes_\RR\CC=S\oplus\Sq=S\oplus S^*=T^*_S$
where $T^*$ denotes the complex algebraic cotangent bundle. With
$G:=K_\CC$ consider the complex algebraic $G$\_variety
$X:=G\times^{H_\CC}S$. Then $T^*_X=G\times^{H_\CC}Z$ where $Z$ is the
restriction of $T^*_X$ to the fiber $S\subseteq X$. Let
$\fq\subseteq\fg$ be a $H_\CC$\_stable complement. Then for every
point $s\in S$ we have a splitting $T_{X,s}=\fq s\oplus
T_{S,s}$. Since $\fq^*\cong\fh^\perp_\CC$ we obtain $Z=\fq^*\times
T^*_S\cong(\fh^\perp\times S)_\CC$. Therefore, we have constructed an
isomorphism $T^*_X\pf\sim M_\CC$. It is easy to see that the moment
map on $M$ induces just the natural cotangent bundle moment map on
$T^*_X$. Thus we have a commutative diagram
$$
\matrix{M&\into&T^*_X\cr
\downarrow\Rechts\mu&&\downarrow\Rechts{\mu_\CC}\cr
\fk^*&\into&\fg^*\cr}
$$
In \cite{CM}~3.3, a subspace $\fa'_\CC\subseteq\ft^*_\CC$ was
constructed such that the closure of the image of $\mu_\CC$ equals
$\overline{G\fa'_\CC}$. The complexified space $\fa^0_\CC$ has the
same property (see \cite{fa0}). Moreover, both have the same image in
$\ft^*/W$ thus they are conjugate by an element $w$ of $W$. The
construction of $\fa'_\CC$ in \cite{CM} depended on the choice of a
Borel subgroup $B$. Thus replacing $B$ by $wBw^{-1}$ we may assume
$\fa'_\CC=\fa^0_\CC$.

Now consider the fiber product $\tilde
T=T^*_X\times_{\ft^*/W}\fa^0_\CC$. Then the group $W(\fa^0)$ acts
transitively on the irreducible components of $\tilde T$ (or at least
on those which map dominantly to the image of $T^*_X$ in
$\ft^*_\CC/W$). There is one such component which is distinguished.
Consider the fiber product $\tilde M:=M\times_{\ft^*/W}\ft^*_+$ which
is an irreducible semi\_algebraic space. Thus the Zariski closure of
the image of $\tilde M\pfeil\tilde T$ is an irreducible component
$\hat T$. Observe that $\hat T$, being the closure of a set of real
points, is defined over $\RR$. The component $\hat T_X$ in
\cite{CM} depended again on the choice of $B$. So changing it we may
assume that $\hat T=\hat T_X$.

Now we let $W(x)$ be the stabilizer of $\hat T$ in $W(\fa^0)$. It is
generated by reflections by \cite{WuM}~6.6. The morphism $\hat
T/W(x)\pfeil T^*_X$ is finite and birational, hence an
isomorphism. Since $\hat T$ maps to $\fa^0_\CC$ we obtain morphisms
$T^*_X\pfeil\fa^0_\CC/W(x)$ and
$\Phi_\CC:T^*_X\pfeil\tilde{G\fa^0_\CC}/W(x)$ where
$\tilde{G\fa^0_\CC}$ is the normalization of the closure of
$G\fa^0_\CC$ in $\fg^*\times\fa^0_\CC$.

By definition, the variety $M_X$ of \cite{WuM} is normal and finite
over $\overline{G\fa^0_\CC}/W(x)$. Moreover $M_X\pfeil\fa^0_\CC/W(x)$ is
the categorical quotient. The generic orbits are closed and the
generic stabilizers of $\overline{G\fa^0_\CC}/W(x)$ are connected. We
conclude that $M_X=\tilde{G\fa^0_\CC}/W(x)$.

By \cite{HC}~9.4, $\Phi_\CC$ induces an isomorphism
$$
\CC[\tilde{G\fa^0_\CC}/W(x)]\pf\sim\|Col|T^*_X.
$$
Since all maps and varieties are defined over $\RR$ we get an
algebraic map between real points
$\Phi:M\pfeil\tilde{K\fa^0}/W(x)$ which induces an isomorphism
$\RR[\tilde{K\fa^0}/W(x)]\pf\sim\|Col|_{\rm alg}(M)$.

Finally, to get an isomorphism between completions, we use that that
$\RR^{>0}$ acts on $M=G\times^H(\fh^\perp\times S)$ by
$t\cdot[g,(\xi,s)]:=[g,(t^2\xi,ts)]$ which equips all algebras with a
grading. Moreover, $\Phi$ is homogeneous. Since $Gx$ is exactly the
fixed point set of this action, we see that completion at $Gx$ is the
same as replacing the direct sum of homogeneous components by their
direct product. This shows the claim.\qed

We proceed with the following comparison result:

\Lemma Diagramm. Let $M$ and $x$ be as above and let $U$ be some $K$\_stable
open neighborhood of $x$. Then there is an isomorphism
$\|Col|_{\rm alg}(M)^\wedge_{Kx}\pf\sim\|Col|(U)^\wedge_{Kx}$.

\Proof: The image of $\RR[M]^K\into\Cinf(U)^K$ is dense
(Stone\_Weierstra{\ss}), thus we obtain $\|Col|_{\rm
alg}(M)\into\|Col|(U)$ and $\|Col|_{\rm
alg}(M)^\wedge_{Kx}\into\|Col|(U)^\wedge_{Kx}$.

To show surjectivity we again use the $\RR^{>0}$\_action on $M$. For
any subspace $S$ of $\Cinf(M)$ let $S_d$ be its subset of elements
which are homogeneous of degree $d$ with respect to this $\RR^{>0}$
action. Let $I\subseteq\Cinf(U)$ be the ideal of functions
vanishing in $K/H$. Then one sees easily
$\Cinf(U)/I^{d+1}=\oplus_{i=0}^d\Cinf(M)_d$. Thus
$\Cinf(U)^\wedge_{Kx}=\prod_{d=0}^\infty\Cinf(M)_d$. The Poisson bracket on
$M$ is homogeneous of degree $-2$, i.e.,
$\{\Cinf(M)_a,\Cinf(M)_b\}\subseteq\Cinf(M)_{a+b-2}$. This implies that
an element of $\Cinf(U)^\wedge_{Kx}$ is collective if and only if all its
homogeneous components are. An analogous result holds for
$\RR[M]^\wedge_{Kx}$. Thus, it suffices to show that $\|Col|_{\rm
alg}(M)_d\pfeil\|Col|(M)_d$ is surjective for every $d\in\NN$.

For an irreducible representation $\eta$ of $K$ and a $K$\_module $V$
let $V^\eta$ denote its $\eta$\_isotypic component. Let $\cE\pfeil
G/H$ be a $K$\_equivariant vector bundle of finite rank. Then
$\Gamma(G/H,\cE)^\eta$ is finite\_dimensional where $\Gamma$ denotes
differentiable sections. This implies in particular
$\|dim|\Cinf(M)_d^\eta<\infty$. From Stone\_Weierstra{\ss} we obtain
$\RR[M]_d^\eta\pf\sim\Cinf(M)_d^\eta$. Thus also $\|Col|_{\rm
alg}(M)_d^\eta\pf\sim\|Col|(M)_d^\eta$.

Since $\RR[\tilde{K\fa^0}/W(x)]$ is a finite $\RR[\fk^*]$\_module we
have $\|dim|\RR[\tilde{K\fa^0}/W(x)]|_d<\infty$. Thus, by
\cite{Comparison}, also $\|Col|_{\rm alg}(M)_d$ is finite
dimensional. We conclude that $\|Col|_{\rm
alg}(M)_d\pf\sim\oplus_\eta\|Col|(M)_d^\eta=:\|Col|(M)_d^{\rm fin}$ is
finite dimensional. But $\|Col|(M)_d^{\rm fin}$ is dense in
$\|Col|(M)_d$ (Peter\_Weyl). Thus equality holds and we obtain
$\|Col|_{\rm alg}(M)_d\pf\sim\|Col|(M)_d$.
\qed

From this we deduce an analogue of \cite{Main} on the level of formal
power series.

\Lemma Formell. Let $M$ be any Hamiltonian $K$\_manifold. For every
$x\in M$ there is a (unique) subgroup $W(x)\subseteq
W(\fa^0)_{\psi(x)}$ which is generated by reflections and an open
$K$\_stable neighborhood $U$ of $x$ such that:
\item{\ri} The map $\Phi=\tilde\mu/W(x):M\pfeil\tilde{K\fa^0}/W(x)$ is
differentiable on $U$ and induces an isomorphism
$\Cinf(\tilde{K\fa^0}/W(x))^\wedge_{K\Phi(x)}\pf\sim\|Col|(U)^\wedge_{Kx}$.
\item{\rii} $U$ carries an analytic structure such that
$\Phi|_U$ is analytic.\Par

\Proof: We may change $x$ in its orbit such that $\mu(x)\in\ft^*_+$
Let $L:=K_{\mu(x)}\subseteq K$. Then $Lx$ is isotropic in $M(L)$ (see
\cite{Cross}). Thus \cite{Comparison} provides us with a group
$W(x)\subseteq W_L$ such that $\Phi_L:M(L)\pfeil
Y:=\tilde{L\fa^0}/W(x)$ is differentiable and even analytic in a neighborhood
of $Lx$. We obtain a commutative diagram
$$
\matrix{
K\times^LM(L)&\pfeil&M\cr
\downarrow\Rechts{\Phi'}&&\downarrow\Rechts\Phi\cr
K\times^LY&\pfeil&\tilde{K\fa^0}/W(x)\cr}
$$
where the horizontal arrows are diffeomorphisms (\cite{Cross},
\cite{Scheibe1}) and $\Phi'=K\times^L\Phi_L$ is differentiable in a neighborhood
of $Kx$ or its image. Thus also $\Phi$ is differentiable and analytic in
a neighborhood of $Kx$.

By \cite{Comparison} and \cite{Diagramm} we have an isomorphism
$\Cinf(Y)^\wedge_{\Phi'(x)}\pf\sim\|Col|(M(L))^\wedge_{Lx}$.  We need
a parameter dependent version of it. Let $\|Col|_K(M(L))$ be the set
of differentiable functions $h_k(x)=h(k,x)$ on $K\times M(L)$ which are
collective for every fixed $k$. Choose functions $f_i\in\Cinf(Y)$
whose Taylor series form a topological basis of
$\Cinf(Y)^\wedge_{\Phi'(x)}$. Then the image of $h_k$ in
$\Cinf(M(L))^\wedge_{Lx}$ can be written as
$\sum_ia_i(k)f_i(\Phi(x))$. The coefficients $a_i$ are unique and can
in fact be expressed in terms of finitely many derivatives of
$h_k$. Thus they are differentiable functions in $k$. This shows
$$
\Cinf(K\times Y)^\wedge_{K\times\Phi'(x)}\pf\sim
\|Col|_K(M(L))^\wedge_{K\times Lx}
$$
and therefore
$$
\Cinf(K\times^L Y)^\wedge_{K\Phi'(x)}\pf\sim
\|Col|_K(M(L))^{\wedge\,L}_{K\times Lx}.
$$
A differentiable function $h$ on $K\times^LM(L)$ is collective if and only if
it is locally constant on fibers of $\Phi'$
(\cite{KarLer}). Therefore, for every $k\in K$ the restriction $h_k$ of $h$
to $[k,M(L)]\cong M(L)$ is collective. Thus we get an inclusion
$$
\|Col|(K\times^LM(L))^\wedge_{Kx}\into
\|Col|_K(M(L))^{\wedge\,L}_{K\times Lx}.
$$
This implies
$$
\Cinf(K\times^L Y)^\wedge_{K\Phi'(x)}\pf\sim
\|Col|(K\times^LM(L))^\wedge_{Kx}
$$
which finishes the proof.\qed

Now we come to a key lemma which compares the little Weyl groups
$W(x)$ of nearby points. The proof rests on an algebraic statement
which is proved in section \cite{Geometry II}.

\Lemma Coherence. Let $M$ be a Hamiltonian $K$\_manifold. Then every $x\in
M$ has an open $K$\_stable neighborhood $U$ such that
$W(y)=W(x)_{\psi(y)}$ for all $y\in U$.

\Proof: With the help of the the Cross Section \cite{Cross} and the
Slice \cite{Slice} we can assume we are in the situation
$M=K\times^H(\fh^\perp\times S)$ and $x=[1,(0,0)]$. Then we have an
algebraic morphism $\phi:M\pfeil\fa^0/W(x)$ such that the formal
power series ring at $0$ in $\fa^0/W(X)$ identifies with the formal
invariant collective functions along $Kx$.

For $y\in M$ and $u:=\psi(y)$, put $W_y:=W(x)_u$. Let $u'$ be the
image of $u$ in $\fa^0/W_y$. Then $\fa^0/W_y\pfeil
\fa^0/W(x)$ is invertible in a neighborhood of $u'$, i.e.,
$M\pfeil\fa^0/W_y$ is differentiable in $Ky$. As in \cite{Minimal} this implies
$W(y)\subseteq W_y$.

The maps
$M\pfeil\fa^0/W(y)\pfeil\fa^0/W_y\pfeil\fa^0/W(x)$ are all differentiable in
$y$. So they induce
$$
\RR[\fa^0]^{\wedge W(x)}_u\pf\sim\RR[\fa^0]^{\wedge W_y}_u
\into\RR[\fa^0]^{\wedge W(y)}_u\into\RR[M]^\wedge_{Ky}.
$$
It suffices to show that the second homomorphism is an isomorphism.
For this, we complexify. Then $Gy$ is a closed orbit in $T_X^*$ and we
get
$$
\CC[\fa^0]^{\wedge W(x)}_u\pf\sim\CC[\fa^0]^{\wedge W_y}_u
\pfeil\CC[\fa^0]^{\wedge W(y)}_u\into\CC[T_X^*]^\wedge_{Gy}.
$$
Recall, that both $W_y$ and $W(y)$ are generated by reflections. Thus,
if $W(y)\ne W_y$ then the map $\fa^0/W(y)\pfeil\fa^0/W_y$ is ramified
over $u'$. That means that $\fa^0/W_y$ contains a divisor whose
preimage in $\fa^0/W(y)$ is not reduced. Thus the same happens for
the spectra of completions. Hence also $\CC[\fa^0]^{\wedge
W(x)}_u\pfeil\CC[T_X^*]^\wedge_{Gy}$ has non\_reduced fibers in
codimension one. Since the completion of a reduced algebra of finite
type is reduced we obtain a contradiction to \cite{Reduced}. Thus
$W_y=W(y)$.\qed

\Lemma Weyl. Assume $M$ to be convex. Let $W_M\subseteq W(\fa^0)$ be
the group generated by all $W(x)$, $x\in M$. Then $W(x)=(W_M)_{\psi(x)}$
for all $x\in M$.

\Proof: Let $S:=\psi(M)$. Since $M$ is convex, $\psi$ has connected
fibers. Thus \cite{Coherence} implies that $W(u):=W(x)$ (with
$u=\psi(x)$) is well defined. Since $W_M$ is generated by reflections
and since $S\subseteq\ft^*_+$ there is a unique Weyl chamber $\fa^0_+$
with respect to $W_M$ which contains $S$. Let $\Delta$ be a root
system attached to $W_M$ (we consider roots as affine functions on
$\fa^0$).

Also every $W(u)$ is a reflection group giving rise to a root
subsystem $\Delta(u)\subseteq\Delta$.  The chamber $\fa^0_+$
determines a chamber for $W(u)$ and therefore a set
$\Sigma(u)\subseteq\Delta(u)$ of simple roots. Let $w\in S$ and $w'\in
S$ any nearby point. Since $w'$ is in the fundamental chamber for
$W_M$ it is also in the fundamental chamber for $W(w)$. Thus, by
\cite{Coherence}, $\Sigma(w')=\{\gamma\in\Sigma(w)\mid\gamma(w')=0\}$.

Let $\Sigma$ be the union of all $\Sigma(u)$, $u\in S$. It suffices to
show that $\Sigma$ is the set of simple roots for $W_M$ because then
the isotropy groups in $\fa^0_+$ are generated by simple
reflections. Since $\Delta$ is certainly generated by the reflection
corresponding to $\Sigma$ it suffices to show that $\Sigma$ is a set
of simple roots for some root system. This is equivalent to
$(\alpha,\beta)\le0$ for all $\alpha\ne\beta\in\Sigma$ where
$(\cdot,\cdot)$ is some $W$\_invariant scalar product on $\Delta$. So
assume $\alpha\in\Sigma(u)$, $\beta\in\Sigma(v)$. Then we have
$\alpha(u)=\beta(v)=0$. Moreover, $u,v\in\fa^0_+$ implies
$\alpha(v)\ge0$, $\beta(u)\ge0$.

Assume first $\alpha(v)=0$. Since $M$ is convex, the line segment $[u,v]$ is
entirely contained in $S$. Moreover, $\alpha$ vanishes on it. This
implies that $\alpha\in\Sigma(w)$ for all $w\in[u,v]$. In particular,
$\alpha,\beta\in\Sigma(v)$ which shows $(\alpha,\beta)\le0$.

The case $\beta(u)=0$ is handled the same way. Thus we are left with
the case where both $\alpha(v)$ and $\beta(u)$ are positive. Let
$\gamma:=\alpha-\beta$. Then $\gamma(u)<0$ and $\gamma(v)>0$. Hence
$\gamma$ is not a root which implies $(\alpha,\beta)\le0$.\qed

\noindent{\it Proof of \cite{Main}:} The group $W_M$ of \cite{Weyl}
certainly satisfies \rii{}. We have $W(x)\subseteq W_M$ for all $x\in
M$. Thus $\Phi=\tilde\mu/W_M$ is differentiable by \cite{Formell}\ri. It
remains to show that every $h\in\|Col|(M)$ is pull\_back from
$\Phi(M)$. For this, we verify the conditions of \cite{FormPush}. Put
$Y:=\Phi(M)$. Since $Y=\nu^{-1}(K\psi(M))/W_M$ it is locally
semi\_algebraic (\cite{ConvProp}\rii). Thus condition \ri{} is
satisfied. Condition \rii{} follows from the convexity of $M$ and
$\tilde\mu(M)/K=\psi(M)$ (\cite{Quotient}). Condition
\riii{} is \cite{Formell}\rii. Finally, condition \riv{} follows again
from $\tilde\mu(M)/K=\psi(M)$ and from the convexity of $M$.

Thus it remains to show that $h$ pushes down formally. For every $x\in
M$ we have $W(x)=(W_M)_{\psi(x)}=(W_M)_{\tilde\mu(x)}$
(\cite{Weyl}). Thus
$$
\Cinf(\tilde{K\fa^0}/W_M)^\wedge_{K\Phi(x)}\pf\sim
\Cinf(\tilde{K\fa^0}/W(x))^\wedge_{K\Phi(x)}.
$$
The assertion follows from \cite{Formell}\ri{}.\qed

\beginsection Alg Ham. Hamiltonian actions in the algebraic category

In the proof of \cite{Coherence} we crucially used a result on complex
algebraic cotangent bundles. The purpose of the next sections is
to provide a proof of this result. All varieties, groups and maps will
be complex algebraic. The group $G$ will always be connected and
reductive. First we develop algebraic versions of the two structure
theorems for Hamiltonian group actions in the algebraic
setting. First, we consider the equivariant Darboux\_Weinstein
Theorem:

\Theorem Darboux-Weinstein. For $i\in\{1,2\}$ let $Z_i$ be a smooth,
affine, symplectic $G$\_variety and $Y_i\subseteq Z_i$ a closed
$G$\_orbit. Let $\phi:Y_1\pfeil Y_2$ be a $G$\_isomorphism and suppose
that the induced isomorphism $d\phi:T_{Y_1}\pf\sim\phi^*T_{Y_2}$
between tangent bundles extends to a $G$\_isomorphism of symplectic
vector bundles $T(\phi):T_{Z_1}|_{Y_1}\pfeil
\phi^*T_{Z_2}|_{Y_2}$. Let $\hat Z_i$ be the spectrum of the
completion of $\CC[Z_i]$ along $Y_i$. Then:
\item{\ri}There exists an isomorphism $\hat\phi:\hat Z_1\pfeil\hat Z_2$ of
symplectic $G$\_schemes which induces $\phi$ and $T(\phi)$.
\item{\rii}Let $\mu_i:Z_i\pfeil\fg^*$ be moment maps with
$\mu_1|_{Y_1}=\phi^*\mu_2|_{Y_2}$. Let $\hat\mu_i$ be the induced
moment map on $\hat Z_i$. Then $\hat\mu_1=\hat\phi^*\hat\mu_2$.\Par

\Proof: Let $x_1\in Y_1$ and $x_2=\phi(x_1)$. Then the isotropy group
$H=G_{x_1}=G_{x_2}$ is reductive. Hence there are $H$\_stable smooth
subvarieties $S_i$ of $Z_i$ which meet $Y_i$ transversally in
$x_i$. Moreover, we can choose the slices $S_i$ in such a way that
$T(\phi)$ maps $V_1:=T_{x_1}(S_1)$ isomorphically to
$V_2:=T_{x_2}(S_2)$. Thus if we put $Z:=G\times^HV_1$, $Y=G/H\subseteq
Z$, and $\hat Z$ the completion of $Z$ along $Y$ we obtain
$G$\_isomorphisms $\hat Z\pfeil\hat Z_i$ such that $\phi$ and
$T(\phi)$ correspond to the identity on $Y$ and $T_{\hat Z}|Y$,
respectively. Moreover, we get two symplectic forms $\omega_i$ on
$\hat Z$. Our assumptions imply $\omega_1|_Y=\omega_2|_Y$. Under these
conditions we have to show that there is a $G$\_automorphism of $\hat
Z$ which maps $\omega_1$ to $\omega_2$ and which is the identity on
$Y$ and $T_{\hat Z}|_Y$.

The proof of this is completely analogous to the proof of the
classical Darboux\_Weinstein theorem given in \cite{GS}, we just have
to make sure that all constructions work in the algebraic category.

Since $Z$ is a vector bundle it carries a scalar multiplication
$\phi:\A^1\times Z\pfeil Z:(t,x)\mapsto\phi_t(x)$. Since
$\phi(\A^1\times Y)=Y$ this induces a morphism $\hat\phi:\hat
Z_t:=(\A^1\times Z)^\wedge_{\A^1{\times}Y}\pfeil \hat Z$. The regular
functions on $\hat Z_t$ can be described as follows. Let $y_1,\ldots,
y_s$ be generators of the sheaf of sections of $G\times^HV_1^*$. Then

\smallskip
\item{$(*)$}{\it every function on $\hat Z_t$ is a formal power
series in the $y_i$ with coefficients in $\CC[Y][t]$.}

\smallskip\noindent
The scalar multiplication also induces an Euler vector field $\xi$ on
$Z$, $\hat Z$, and $\hat Z_t$. Let $\sigma:=\omega_1-\omega_0$ which
is a closed $2$\_form on $\hat Z$. Then we can construct the $1$\_form
$$
\beta:=\int_0^1\iota_\xi(\phi^*\sigma) dt.
$$
By $(*)$, the integration is over polynomials in $t$, hence
well\_defined. Moreover, since both $\xi$ and $\sigma$ vanish in $Y$,
the form $\beta$ vanishes quadratically.

Next we form
$\omega_t:=\omega_0+t\sigma=(1-t)\omega_0+t\omega_1$. Since $\omega_0$
is non\_degenerate and $\sigma|_Y=0$, this is a non\_degenerate closed
$2$\_form on $\hat Z_t$ (relative to $\A^1$). Thus we can define a
vector field $\eta_t$ on $\hat Z_t$ by $\eta_t\haken\omega_t=-\beta$.

I claim that one can integrate the time dependent vector field
$\eta_t$ on $\hat Z$. For this, it suffices to show that the
derivation ${\partial\over\partial t}+\eta_t$ acts pointwise
topologically nilpotently on $\CC[\hat Z_t]$. But $\eta_t$ also
vanishes quadratically in $Y$. Thus it strictly increases the order of
vanishing along $Y$. On the other hand, a sufficiently high power of
${\partial\over\partial t}$ also increases the order of vanishing by
$(*)$.

Thus, we obtain a one\_parameter family $f_t$ of morphisms of $\hat Z$
into itself such that $f_0=\|id|$. Moreover, $f_t$ is the identity on
$Y$ and $T_{\hat Z}|_Y$ since $\eta_t$ vanishes quadratically along $Y$. 
This implies that $f_t$ is an automorphism for every $t$, in
particular, for $t=1$. Furthermore, everything being $G$\_equivariant,
also $f_1$ commutes with the $G$\_action.

Finally, the point of all this is of course the relation
$\omega_0=f_1^*\omega_1$. But this is a formal consequence of the
identity (22.1) in \cite{GS} and several other identities which hold
also in the algebraic setting. One way to see this is: they hold
for holomorphic forms and those are dense in the forms with
coefficients in formal power series.

This shows {\ri}. Part {\rii} follows from the fact that a moment map is
unique up to a translation in $\fg^*$.\qed

\noindent As consequence we obtain the algebraic analogue of
\cite{Slice}.

\Corollary SliceAlg. Let $Gx\subseteq Z$ be a closed isotropic orbit. Put
$S_x:=(\fg x)^\perp/\fg x$ which is a symplectic representation of
$G_x$. Let $u_x:=\mu(x)\in(\fg^*)^G$. Then the triple $(G_x,S_x,u_x)$
determines a formal neighborhood of $Gx$ uniquely up to Hamiltonian
isomorphism. Conversely, any such triple occurs.

We temporarily identify $\fg^*$ with $\fg$ by means of a
$G$\_invariant scalar product. Then every element $v\in\fg$ has a
Jordan decomposition $v=v_s+v_u$.  Let $L\subseteq G$ be a Levi
subgroup, i.e., the centralizer of a semisimple element in $\fg$. Let
$W_L$ be the Weyl group of $L$. Consider $\ft^r:=\{v\in\ft^*\mid
W_v\subseteq W_L\}$. In other words $\ft^r$ is obtained from $\ft$ by
removing all reflecting hyperplanes which don't belong to $W_L$. Hence
it is an open subset of $\ft^*$. It has the property that
$G_v\subseteq L$ whenever $v\in L\ft^r$.

\Lemma Scheibe. The morphism
$$
\lambda:\quad G\times^L(\fl\times_{\ft/W_L}\ft)\pfeil\fg\times_{\ft/W}\ft
:\quad[g,(v,w)]\mapsto(gv,w)
$$
is an isomorphism over $\ft^r$. Moreover, if $\fl^r\subseteq\fl$
is the preimage of $\ft^r/W_L\subseteq\fl\mod L$ then we have
$$
G\times^L\fl^r\pf\sim\fg\times_{\ft/W}\ft^r/W_L
$$

\Proof: It is well known (see~?) that both sides are normal
varieties. So, by Zariski's main theorem (or just the Richardson
lemma), it suffices to show that $\lambda$ is bijective over $\ft^r$.

Assume $\lambda([g,(v,w)])=\lambda([g',(v',w)])$. Since $v$ and $v'$
have the same image as $w$ in $\ft/W=\fl\mod L$ there is $l\in L$ such
that $v_s=lv'_s\in Lw$. Furthermore, $gv=g'v'$ implies
$gv_s=g'v'_s$. Thus $lg^{\prime-1}g\in G_{v_s}$. Since $v_s\in L\ft^r$
we have $g^{\prime-1}g\in L$. Thus
$[g,(v,w)]=[g',(g^{\prime-1}gv,w)]=[g',(v',w)]$.

Conversely, let $v\in\fg$ and $w\in\ft^r$ have the same image in
$\ft/W$. Then there is $g\in G$ with $v_s=gw$. Hence
$u:=g^{-1}v\in\fg_w\subseteq\fl$ and we obtain
$\lambda([g,(u,w)])=(v,w)$.

The second isomorphism follows from the first by dividing out $W_L$ on
both sides.
\qed

Next, we consider the cross section theorem.

\Theorem CrossAlg. Let $Z$ be an algebraic Hamiltonian $G$\_manifold with
moment map $\mu$. For
a Levi subgroup $L\subseteq G$ put $Z(L):=\mu^{-1}(\fl^r)$ and
$Z_L:=G\times^LZ(L)$. Assume, $Z(L)$ is not empty. Then
\item{\ri} The set $Z(L)$ is a Hamiltonian $L$\_variety (possibly
disconnected): its symplectic form
is the restriction of that of $Z$; its moment map is
$Z(L)\pfeil\fl^r\into\fl^*$.
\item{\rii} The canonical morphism $Z_L\pfeil Z$ is an {\'e}tale Poisson
morphism.
\item{\riii}If $Z$ is affine then $Z(L)$ and $Z_L$ are affine as
well. Moreover, let $Y\subseteq Z_L$ be a closed orbit. Then its image
$Y'$ in $Z$ is also closed and $Y\pfeil Y'$ is an isomorphism.\Par

\Proof: From \cite{Scheibe} we obtain an isomorphism
$$
Z_L=G\times^LZ(L)\pf\sim Z\times_{\ft/W}\ft^r/W_L.
$$
Then $Z_L\pfeil Z$ is {\'e}tale since $\ft^r/W_L\pfeil\ft^*/W$ is. The
other assertions are either trivial or can be deduced exactly is in
the classical case (\cite{Cross}).\qed

\Remark: The Slice \cite{SliceAlg} works best for closed orbits $Gx$
with $\mu(x)=0$. On the other hand, \cite{CrossAlg} allows reduction
to the case where $\mu(x)$ is nilpotent (take $L$ to be the
centralizer of $\mu(x)_s$). Thus, unlike in the classical case, there
is a gap between the ranges of applicability of these theorems.

\beginsection Cotangent. The moment map of complex cotangent bundles:
connectedness of fibers

Again, $G$ is a connected reductive group, $B\subseteq G$ a Borel
subgroup, $U\subseteq B$ its unipotent radical, and $T\subseteq B$ a
maximal torus. Let $X$ be a smooth $G$\_variety. Then the cotangent
bundle $T_X^*$ is equipped with a canonical symplectic structure and a
moment map $\mu:T_X^*\pfeil\fg^*$. Then, in \cite{WuM}, we have
constructed a subspace $\fa^*$ of $\ft^*$, a subquotient $W_X$ of the
Weyl group $W$ and a morphism $\phi:T^*_X\pfeil\fa^*/W_X$ which
factors $T^*_X\pfeil\fg^*\mod G=\ft^*/W$. It is known to be
flat. Almost by definition its generic fibers are connected. We prove
a refinement.

\Theorem ConnGen. Let $X$ be a smooth $G$\_variety. Then all fibers of
$\phi:T_X^*\pfeil\fa^*/W_X$ are connected.

\Proof: Let us assume first that $X=G/H$ is homogeneous.
Then there is a nice way to compactify the moment map $\mu$. The Lie
algebra $\fh$ of $H$ corresponds to an $H$\_fixed point in the
Gra{\ss}mannian $\|Gr|(\fg)$ of $\fg$. Thus we get a $G$\_morphism
$X\pfeil\|Gr|(\fg)$. Choose any $G$\_equivariant compactification
$X\into\XS$ of $X$. Then $G/H$ maps diagonally into
$\|Gr|(\fg)\times\XS$. Let $\Xq$ be the normalization of closure of
the image. Thus, by construction, $\Xq$ is a compactification of $X$
such that the canonical map $X\pfeil\|Gr|(\fg)$ extends to
$\Xq\pfeil\|Gr|(\fg)$.

Let $V\pfeil\|Gr|(\fg)$ be the vector bundle whose fiber over a point
$\fm\subseteq\fg$ is $\fm^\perp\subseteq\fg^*$. Let
$\Tq:=\Xq\times_{\|Gr|(\fg)}V$. Then, by construction, the restriction
of $\Tq$ to $X$ is $T_X^*$. Since $V\subseteq\|Gr|(\fg)\times\fg^*$ we
obtain a morphism $\overline\mu:\Tq\pfeil\fg^*$ which extends
the moment map. Moreover, by construction, $\overline\mu$ is proper.

Let $\Tq\pfeil M\pfeil\fg^*$ be the Stein factorization of
$\overline\mu$. Then $M$ is an affine $G$\_variety with
$M/\!\!/G=\fa^*/W_X$. In particular, the fibers of $M_X\pfeil\fa^*/W_X$ are
connected. By Zariski's connectedness theorem, all fibers of
$\Tq\pfeil M$ are connected. This implies that also all fibers of
$\overline\phi:\Tq\pfeil\fa^*/W_X$ are connected.

The fiber $\Tq_x$ of $\Tq$ over $x\in\Xq$ is
$\fm(x)^\perp\subseteq\fg^*$ where $\fm(x)\subseteq\fg$ is a limit of
conjugates of $\fh$. It follows that $\fm$ is also an algebraic Lie
algebra with the same rank and complexity as $\fh$ (see
\cite{WuM}~2.5). This implies that the map $\fm(x)\pfeil\fa^*/W_X$ is
equidimensional (\cite{WuM}~6.6). Thus also
$\overline\phi:\Tq\pfeil\fa^*/W_X$ is equidimensional, hence flat,
since both domain and range are smooth varieties.

Consider now the fiber $F:=\overline\phi^{-1}(u)$ over a point
$u\in\fa^*/W_X$. Since all maps $\fm(x)^\perp\pfeil\fa^*/W_X$ are
equidimensional, no irreducible component of $F$ is contained in
$\partial\Tq:=\Tq\setminus T^*_X$. Let $C_1$, $C_2$ be irreducible
components of $\phi^{-1}(u)$. Then their closures $\Cq_i$ are
irreducible components of $F$. Since $F$ is locally a complete
intersection it is connected in codimension one (\cite{Ha}). Thus, to
show that $\phi^{-1}(u)$ is connected it suffices to prove the
following claim: assume $\Cq_1\cap\Cq_2$ is non\_empty, of codimension
one in $\overline\phi^{-1}(u)$, and contained in $\partial\Tq$. Then
$C_1\cap C_2\ne\emptyset$.

Consider the subset $\Xq_0\subseteq\Xq$ of $x\in X$ where
$\|dim|(\fm(x)+\fu)$ and $\|dim|(\fm(x)+\fb)$ are maximal. Since
$\fm(x)$ is an algebraic Lie algebra of the same complexity and rank
as $\fh$, the set $X_0$ meets all $G$\_orbits of $\Xq$. Now consider
the subbundels $\cB\subseteq\cU\subseteq\Tq|_{\Xq_0}$ whose fibers at
$x$ are $\cB_x=(\fm(x)+\fb)^\perp$ and
$\cU_x:=(\fm(x)+\fu)^\perp$. The quotient $\cU/\cB$ can be identified
with the trivial bundle $\fa^*\times\Xq_0$. Moreover, the restriction
of $\overline\phi$ to $\cU$ is just the projection onto $\fa^*$
followed by the quotient map $\fa^*\pfeil\fa^*/W_X$ (see
\cite{WuM}~6.2).  Thus, the intersection $\Cq_i\cap\cU$ is contained
in the preimage of a finite subset $S\subseteq\fa^*$. On the other
hand, its codimension is less or equal $\|dim|\fa^*$. We conclude that
$\Cq_i\cap\cU$ {\it equals\/} the preimage of some $S_i$. Therefore,
we are done if we can show that $\Cq_1\cap\Cq_2\cap\cU$ is not empty
since then $S_1\cap S_2\ne\emptyset$.

Since $I:=\Cq_1\cap\Cq_2\subseteq\partial\Tq$ and since it is of
codimension one in $\overline\phi^{-1}(u)$, it must be the union of irreducible
components of $\partial\Tq\cap\overline\phi^{-1}(u)$. Thus there is
$x\in\Xq_0$ such that $I\cap\fm(x)^\perp$ is a union of components of
$\phi_x^{-1}(u)$ where $\phi_x$ is the restriction of $\overline\phi$
to $\fm(x)^\perp=\Tq_x$. Write $\cU_x=\cB_x\oplus\fa^*(x)$. Then we
are done with the following lemma:

\Lemma XYZ. Every irreducible component $C$ of $\phi_x^{-1}(u)$ meets
$\fa^*(x)$.

\Proof: We use the $\CC^*$\_action by scalar multiplication. Then
$\phi_x^{-1}(\CC^*u)$ is closed in
$\fm(x)^\perp\setminus\phi_x^{-1}(0)$. Thus
$\overline{\CC^*C}\subseteq\CC^*C\cup\phi_x^{-1}(0)$. Since
$\overline{\CC^*C}\cap\fa^*(x)\ne\emptyset$ and
$\|codim|\CC^*C=\|dim|\fa^*-1$ we have
$\|dim|\overline{\CC^*C}\cap\fa^*(x)\ge1$. Moreover, this intersection
is not contained in $\phi_x^{-1}(0)$. Thus
$\CC^*C\cap\fa^*(x)\ne\emptyset$ which implies the claim.\qed

This settles the case where $X$ is homogeneous. In the general case
let $X_0\subseteq X$ be non\_empty, open, $G$\_stable such that the
orbit space $X_1:=X_0/G$ exists. Then the morphism
$T_{X_0}^*\pfeil\fa^*/W_X$ factors through the relative cotangent
bundle $T_{X_0/X_1}^*$. The discussion of the homogeneous case implies
that the fibers of $T_{X_0/X_1}^*\pfeil Y\times\fa^*/W_X$ are
connected. Thus the assertion of the theorem holds for $X_0$.

Let $Y\subseteq X$ be an orbit and let $\pi:T^*_X\pfeil X$ be the
projection. The restriction of $\phi$ to $\pi^{-1}(Y)$ factors through
$T^*_Y$. Let $C$ be any component of
$\pi^{-1}(Y)\cap\phi^{-1}(u)$. Then, by \cite{XYZ}, there is $x\in Y$
and a subspace $\fa^*_Y(x)\subseteq T^*_{X,x}$ such that
$\fa^*_Y(x)\cap C\ne\emptyset$. Now we are applying the local
structure theorem to $Y$ (\cite{IB}~2.10): there is a parabolic
subgroup $P$ with Levi part $L$ and an affine subvariety $Z$ such that
$P_u\times Z\pfeil X:(u,z)\mapsto uz$ is an open embedding, with
$Z\cap Y\ne\emptyset$ and $(L,L)$ acts trivially on $Z\cap
Y$. Moreover $(\|Lie|L/L_x)^*=\fa^*_Y$. For $x\in Z\cap Y$ we can find
a $L_x$\_stable slice $S$ in $Z$. Then for every $z\in S$ choose a
complement $\fa^*_Y(z)$ of $(T_{Z,z})^\perp$ in
$(T_{S,z})^\perp\subseteq T^*_{X,z}$ which can be canonically
identified with $\fa^*_Y$. Moreover, the restriction of $\phi$ to
$\fa^*_Y(z)$ is just the composition
$\fa^*_Y(z)\into\fa^*\pfeil\fa^*/W_X$. This shows that the
intersection $C\cap\fa^*_Y(x)$ can be moved continuously into
$X_0$. Thus, $\phi^{-1}(u)$ is connected.\qed

\beginsection Geometry II. The moment map of complex cotangent bundles:
reducedness of fibers

Next we investigate the local structure of the fibers of
$\phi:T_X^*\pfeil S:=\fa^*/W_X$. In \cite{ARC} I have shown how to
interpret it as a moment map: there is a flat abelian group scheme
$\cA_X/S$ which acts on $T_X^*/S$. Moreover, there is an isomorphism
$T^*_S\pfeil\|Lie|\cA_X$ (in particular, the fibers of $\cA_X$ are
$\|dim|\fa^*$ dimensional) such that the following diagram commutes
(where $Z:=T_X^*$):
$$
\matrix{T^*_S\times_SZ&\pf{\phi^*}&T^*_Z\cr
\Links\sim\untenPf&&\Links\sim\untenPf\Rechts\alpha\cr
(\|Lie|\cA_X)\times_SZ&\pf\beta&T_Z\cr}
$$
Here, $\alpha$ is induced by the symplectic structure on $Z$ and
$\beta$ by the action of $\cA_X$ on $Z$. If $\cA_X$ were a constant
group scheme then the diagram above would be equivalent to the
defining property of the moment map.

The commutative diagram above shows in particular, that $\phi$ is
smooth in a point $x$ if and only if its $\cA_X$\_isotropy group is
finite.

The group scheme $\cA_X/S$ is intimately connected with the quotient
map $\fa^*\pfeil S=\fa^*/W_X$. As any commutative linear group, the
fiber $\cA_{X,u}$ over a point $u\in S$ is the direct product of its
semisimple part $\cA_{X,u}^s$ and its unipotent part
$\cA_{X,u}^u$. Choose $x\in\fa^*$ over $u$ and let $W_{X,x}$ be its
isotropy group. Then $\cA_{X,u}^s$ is an open subgroup of
$A_X^{W_{X,x}}$. This determines also the dimension of
$\cA_{X,u}^u$. In particular, $\cA_{X,u}$ is a torus if and only if
$u$ is not in the branching divisor $\Delta$ of $\fa^*\pfeil
S$. Furthermore, $\|dim|\cA_{X,u}^u=1$ if and only if $|W_{X,x}|=2$ if
and only if $u$ is a smooth point of $\Delta$.

\Example GroupSchEx. \rm Take $X=G$ with the action by left
translations. The cotangent bundle $T^*_G$ is trivial and the
$G$\_action on it becomes $g(h,\lambda)=(gh,\lambda)$. The moment map
is $\mu(g,\lambda)=g\lambda$. The action of $\cA_X$ on $T^*_X$
commutes with the $G$\_action and preserves $G$\_invariants. This
implies that for every $a\in\cA_G$ and $\lambda\in\fg^*$ which have
the same image in $\ft^*/W$, there is $c(a,\lambda)\in G$ such that
$a\cdot(g,\lambda)=(gc(a,\lambda)^{-1},\lambda)$. Since the
$\cA_X$\_action also preserves $\mu$, we have $c(a,\lambda)\in
G_\lambda$. Thus, the $\cA_G$\_action on $T^*_G$ is described by a
homomorphism $\cA_G\times_{\ft^*/W}\fg^*\pfeil G\times\fg^*$ such that
$\cA_{G,\pi(\lambda)}\pfeil G_\lambda$. One can show that this is an
isomorphism whenever $\lambda$ is regular. Moreover, we obtain
homomorphisms between Lie algebras
$$
T^*_{\ft^*/W}\times_{\ft^*/W}\fg^*\cong\|Lie|\cA_G\times_{\ft^*/W}\fg^*
\pfeil\fg\times\fg^*=T^*_{\fg^*}
$$
which is just the derivative of $\pi:\fg^*\pfeil\ft^*/W$.

\Theorem AAA. Let $u\in S\setminus\Delta$. Then the fiber of
$\phi:T^*_X\pfeil S$ over $u$ is reduced.

\Proof: Let $\fa^r\subseteq\fa^*$ be
the preimage of $S\setminus\Delta$. Since $\fa^r\pfeil S$ is {\'e}tale,
$Z:=T_X^*\times_S\fa^r$ is a smooth symplectic variety. The pull\_back
$\cA_X\times_S\fa^r$ is just the trivial group scheme with fiber
$A_X$. Thus $A_X$ acts on $Z$. The moment map of this action is the
projection $Z\pfeil\fa^r\subseteq\fa^*$. Since the assertion of the
lemma is stable under {\'e}tale base\_change, we have reduced it to a
statement about moment maps of tori on affine varieties.

Let $z\in\phi^{-1}(u)$. Since $A:=A_X$ is a torus, there is an open
affine $A$\_stable neighborhood of $z$ in which $Az$ is
closed. Moreover, $Az$ is automatically isotropic. Thus, by
\cite{Darboux-Weinstein}, a formal neighborhood of $Az$ is isomorphic
to $Z':=A\times^H(\fh^\perp\oplus S)$ where $H:=A_z$. The fiber
$\phi^{-1}(u)$ is reduced in $Az$ if and only if its completion is
(see \cite{ZaS} VIII~Thm.~31). Thus we may replace $Z$ by $Z'$. But
then it suffices to look at the moment map for the $H^0$\_action on
the vector space $S$.

Since $S$ is symplectic there are linear coordinates
$(q_1,\ldots,q_n,p_1,\ldots,p_n)$ of $S$ and characters
$\chi_1,\ldots,\chi_n$ of $H^0$ such that $H^0$ acts on $S$ by
$t\cdot(q_i,p_i)=(\chi_i(t)q_i,\chi_i(t)^{-1}p_i)$. Let
$d\chi:=\|Lie|\chi_i\in\fh^*:=(\|Lie|H^0)^*$. Then the moment map is
given by $\mu(q_i,p_i)=\sum_iq_ip_id\chi_i$. Thus, $\mu$ is the
composition of $\pi:S\pfeil\CC^n:(q_i,p_i)\mapsto(q_ip_i)$ and the
linear map $\CC^n\pfeil\fh^*:(t_i)\mapsto\sum_it_id\chi_i$. Any fiber
of the latter is an affine space $V$, hence reduced. Moreover, since
all fibers of $\pi$ are reduced the generic fiber of
$\pi^{-1}(V)\pfeil V$ is reduced. Because $\pi$ is flat, this implies
that $\pi^{-1}(V)$ itself is reduced.\qed

It remains to study the fibers of $\phi$ over $\Delta$. There our
results are rather incomplete. A first subcase is:

\Lemma BBB. Assume $X$ to be affine. Let $Gz\subseteq T^*_X$ be a
closed orbit such that $G\mu(z)\subseteq\fg^*$ is {\rm not}
closed. Assume moreover, that $u:=\phi(z)$ is a smooth point of
$\Delta$. Then $\phi^{-1}(\Delta)$ is smooth in $z$.

\Proof: Since $\Delta$ is smooth in $u$, the unipotent part
of $\cA_{X,u}$ is one\_dimensional. The Lie algebra of it corresponds
via the identification $\|Lie|\cA_{X,u}=T^*_{S,u}$ to
$(T_u\Delta)^\perp$. We conclude that $C$ is smooth in $z$ if and only
if $z$ is not a fixed point for $\cA_{X,u}^u$. This is what we are
going to show.

Let $\mu(z)=:\lambda=\lambda_s+\lambda_u$ be the Jordan decomposition
(observe $\fg^*\cong\fg$) and put $L:=G_{\lambda_s}$. We may assume
$\lambda_s\in\fa^*\subseteq\ft^*$ and thus
$\ft^*\subseteq\fl^*$. Let $u'\in\ft^*/W_L$ and $u''\in\ft^*/W_G$ be the
images of $\lambda$. Let $\cA_G$ and $\cA_L$ be the group schemes
associated to $G$ and $L$ (see \cite{GroupSchEx}). Then we obtain the
following commutative diagram (with $Z:=T^*_X$):
$$
\matrix{
\fl&\cong&T^*_\lambda\fl^*&\leftarrow&T^*_{u'}(\ft^*/W_L)
&\cong&\|Lie|\cA_{L,u'}\cr
\uparrow\Rechts{\alpha}&&\uparrow\Rechts{\alpha'}&&\uparrow\Rechts\sim&
&\uparrow\Rechts\sim\cr
\fg&\cong&T^*_\lambda\fg^*&\buildrel\beta\over\leftarrow&T^*_{u''}(\ft^*/W_G)
&\cong&\|Lie|\cA_{G,u''}\cr
\downarrow&&\downarrow&&\downarrow&&\downarrow\cr
T_zZ&\cong&T^*_zZ&\leftarrow&T^*_u(\fa^*/W_X)
&\cong&\|Lie|\cA_{X,u}\cr}
$$
The inclusion $\fc_\fg(\lambda)\subseteq\fl$ implies
$\fl^\perp\subseteq\fc_\fg(\lambda)^\perp=\fg\lambda$. This shows that
the image of $\beta$ is contained in $(\fg\lambda)^\perp=\fl$. This
implies that, if we replace $\alpha$ by the natural inclusion
$\fl\into\fg$ (and remove $\alpha'$) then the diagram still commutes.

Let $q_0$ be any $L$\_invariant non\_degenerate quadratic form on
$\fl^*$ and $q(\tau):=q_0(\tau-\lambda_s)$ for $\tau\in\fl^*$. Then
the 1-form $dq_0$ can viewed as a map $\fl^*\pfeil\fl$. It is affine
linear and maps $\lambda_s$ to $0$ and $\lambda$ to a non\_zero
nilpotent element $\xi\in\fl$. On the other hand, $q_0$, being
$L$\_invariant, is a pull\_back from $\ft^*/W_L$. This show that $\xi$
is the image of some element of $\eta\in\|Lie|\cA_{L,u'}$. The map
$\|Lie|\cA_{L,u'}\pfeil\fl$ is induced by a homomorphism between
algebraic groups (see \cite{GroupSchEx}). Thus we can choose $\eta$ to
be nilpotent. Hence there is also an element of
$\|Lie|\cA_{X,u}^u$ which is mapped to $\xi z\in T_zZ$.

The orbit $Gz$ is closed by assumption. Hence its isotropy group $G_z$
is reductive. Moreover, $\lambda$ and $\xi$ are $G_z$\_fixed. Since
$\fg_z$ does not contain nilpotent elements in its center we conclude
that $\xi\not\in\fg_z$, i.e., $\xi z\ne0$.\qed

For the next result note that $\phi^{-1}(\Delta)$ is a divisor in
$T^*_X$, i.e., purely of codimension one and each component has a
multiplicity attached to it.

\Lemma CCC. Assume $X$ to be affine. Let $K\subseteq T^*_X$ be
a closed, $G$\_stable subvariety of codimension two. Assume $\phi(K)$
is an irreducible component of $\Delta$. Then the multiplicity of all
components of $\phi^{-1}(\Delta)$ which pass through $K$ is the
same. Moreover, it is either one or two.

\Proof: Since we want to apply the Cross Section and the Slice
Theorem, we have to check that the property to be proved can be
detected in formal neighborhoods. Let $f=0$ be a local equation of
$\phi^{-1}(\Delta)$ in $z$. Since the local ring $\cO_{Z,z}$ is
regular, it is an UFD. Therefore $f$ has a prime factorization
$f=ag_1^{m_1}\ldots g_s^{m_s}$. Here, $a$ is a unit, the $g_i$
correspond to the different components of $\phi^{-1}(\Delta)$ passing
through $z$ and the $m_i$ are their multiplicities. We want to show
$m_1=\ldots=m_s=1$ or $m_1=\ldots=m_s=2$. The completion
$\hat\cO_{Z,z}$ is also regular, thus a UFD. The $g_i$ are in general
not prime in $\hat\cO_{Z,z}$ but one can say two things: each
$g_i$ is square\_free in $\hat\cO_{Z,z}$ (since $\hat\cO_{Z,z}/(g_i)$
has no nilpotents, see \cite{ZaS} VIII~Thm.~31), and: if $i\ne j$
then $g_i$ and $g_j$ are coprime in $\hat\cO_{Z,z}$ (since
$\|dim|\hat\cO_{Z,z}/(g_i,g_j)=\|dim|Z-2$). This implies easily that
our assertion can be checked in the completion alone.

Let $Gz\subseteq K$ be a generic closed orbit. If
$G\mu(z)\subseteq\fg^*$ is not closed then $\phi^{-1}(\Delta)$ is even
smooth in $z$ (\cite{BBB}). Thus we have yet to consider the case
where $\xi:=\mu(z)$ is semisimple.  First, we are applying the Cross
Section \cite{CrossAlg} to $Z:=T^*_X$ and $L=G_\xi$. Then it suffices
to prove the assertion for a formal neighborhood of $Lz$ in $Z(L)$ and
the restricted morphism $Z(L)\pfeil S$.

Next we apply \cite{SliceAlg}. Thus, a formal neighborhood of $Lx$ in
$Z(L)$ is isomorphic to a formal neighborhood of the zero\_section of
$Z':=L\times^{L_x}(\fl_x^\perp\oplus S)$ where $S$ is a symplectic
representation of $L_x$. From
$\fg=\fp_u^+\oplus\fl\oplus\fp_u^-$ we obtain isomorphisms of
symplectic $L_x$\_representations:
$$
T_zZ\cong\fp^+_uz\oplus\fp^-_uz\oplus T_xZ(L)\cong
[\fp^+_uz\oplus\fp^-_uz]\oplus[\fl z\oplus (\fl z)^*]\oplus S
$$
It is easy to see that every symplectic representation of a reductive
group is uniquely of the form $V\oplus V^*\oplus\bigoplus_i M_i$ where
$V$ and $V^*$ are isotropic submodules and the $M_i$ are irreducible
and pairwise non\_isomorphic. Since $T_zZ$ has a Lagrangian submodule
we conclude that $S\cong V\oplus V^*$ for some Lagrangian submodules
$V$ and $V^*$. This means that we can identify $Z'$ with the cotangent
bundle $T^*_{X'}$ where $X':=L\times^{L_x}V$.

The morphism $Z'\pfeil S$ factors through $\phi':Z'=T^*_{X'}\pfeil
S':=\fa^*/W_{X'}$. Let $\Delta'$ be the branch divisor of $S'$. If
$\phi'(K)\not\subseteq\Delta'$ then we are in the situation of
\cite{AAA}, i.e., all fibers of $\phi'$ outside $\Delta'$ are
reduced. Moreover all components of the preimage of $\Delta$ in
$S'\setminus\Delta'$ have multiplicity two. Thus, the multiplicities
of all components of $\phi^{-1}(\Delta)$ passing through $K$ are
$2$. (Remark: it will turn out (\cite{Reduced}) that this case
actually never occurs).

We are left with the case $\phi'(K)\subseteq\Delta'$. The group
$W_{X'}$ is contained in the isotropy group of $W_X$ over
$\phi(u)$. Since that point is generic in the codimension one
subvariety $\phi(K)$ we must have $|W_{X'}|=2$. Therefore, the
morphism $S'\pfeil S$ is unramified in the generic point of
$\Delta'$. Thus we may replace the map $Z'\pfeil S$ by $Z'\pfeil
S'$. The moment map on $Z'$ coming from $Z(L)$ differs from the
canonical moment map of $Z'$ as a cotangent bundle by a shift by a
$L$\_fixed element. This means that we can also replace the original
moment map by the cotangent bundle moment map. The conclusion of the
discussion is: we can replace $G, X$ by $L, X'$ and thus assume
that $|W_X|=2$.

Consider the projection $\pi:Z=T_X^*\pfeil X$. I show first that every
component of $\phi^{-1}(\Delta)$ maps dominantly to $X$. Otherwise,
there is a $G$\_stable divisor $Y$ of $X$ such that
$\pi^{-1}(Y)\subseteq\phi^{-1}(\Delta)$. The divisor $Y$ induces a
$G$\_stable valuation $v_Y$ of the function field $\CC(X)$. Since the
restriction of the moment map to $\pi^{-1}(Y)$ factors through the
moment map on $T_Y^*$ we have
$\|dim|\fa^*_Y=\|dim|\Delta=\|dim|\fa^*-1$. This implies that $v_Y$ is
a so\_called central valuation (see
\cite{IB}~7.3). These form a cone in the group
$\cH:=X_*(A_X)$ of one\_parameter subgroups of $A_X$ which is a
$W_X$\_module. The image of $\fa_Y^*\subseteq\fa^*$ in $S$ is
$\Delta$. This implies (\cite{IB}~7.4) that $\sigma\cdot v_Y=-v_Y\in\cH$
where $\sigma$ is the non\_trivial element of $W_X$.

Let $\lambda:\CC^*\pfeil A_X$ be the one\_parameter subgroup attached
to $v_Y$. In \cite{CM}~7.3, I constructed a class of maps
$\psi:A_X\into X$ with the property
$\|lim|_{t\pfeil0}\psi(\lambda(t))\in Y$. With $\psi$ also
$\psi\circ\sigma$ is of this class. Thus
$\sigma(\lambda)=\lambda^{-1}$ implies
$$
\|lim|_{t\pfeil\infty}\psi(\lambda(t))=
\|lim|_{t\pfeil0}(\psi\circ\sigma)(\lambda(t))\in Y
$$
This means that $\psi\circ\lambda:\CC^*\pfeil X$ extends to a
non\_constant morphism from $\P^1$ to $X$ which is absurd since $X$ is
affine.

This shows that we can ``recognize'' all components of
$\phi^{-1}(\Delta)$ by looking at a generic (possibly non\_affine)
orbit $Gy$ of $X$. The restriction of the moment map to $\pi^{-1}(Gy)$
factors through $T^*_{Gy}=G\times^{G_y}\fg_y^\perp$. Now observe that
$\Delta\subseteq S$ is given by the vanishing of a function of degree two.
Since the map $\fg_y^\perp\pfeil S$ is homogeneous, also the intersection of
$\phi^{-1}(\Delta)$ with $\fg_y^\perp$ is given by the vanishing of a {\it
quadratic\/} polynomial. We conclude, that $\phi^{-1}(\Delta)$ is
either irreducible with multiplicity at most two, or it consists of
two components of multiplicity one. That concludes the proof.\qed

\noindent The following is a technical result which was already
used in section~\cite{Proof of main theorem}.

\Theorem Reduced. Let $X$ be a smooth affine $G$\_variety and let
$D\subseteq\fa^*/W_X$ be a prime divisor. Then $\phi^{-1}(D)$ is
reduced.

\Proof: If $D\not\subseteq\Delta$ then $\phi^{-1}(D)$ is reduced by
\cite{AAA}. Thus let $D$ be a component of $\Delta$. Let
$C_1,\ldots,C_s$ be the irreducible components of
$\phi^{-1}(D)$. Define a graph $\Gamma$ whose vertices are the
$C_i$. There is an edge between $C_i$ and $C_j$ if $K_{ij}=C_i\cap C_j$ is
of codimension two in $T^*_X$ and $\phi(K_{ij})=D$. I claim, that $\Gamma$
is connected. Indeed, let $D_0\subseteq D$ be the complement of the
union of all $\phi(K_{ij})\ne D$. Then $\phi^{-1}(D_0)$ has the same
number of components as $\phi^{-1}(D)$ since $\phi$ is flat. Moreover,
$\phi^{-1}(D_0)$ is connected by \cite{ConnGen}. Being locally a
complete intersection one can remove from $\phi^{-1}(D_0)$ all
$K_{ij}$ which are not of codimension two in $T^*_X$ and the result is
still connected (\cite{Ha}). But this means that $\Gamma$ is connected.

It follows now from \cite{CCC} that if $\Phi^{-1}(D)$ is not reduced
then all of its components have multiplicity two. The divisor $D$ is
defined by an equation$f_0=0$. Thus $\phi^{-1}(D)$ is the zero scheme
of $f:=f_0\circ\phi=0$. Thus, if $\phi^{-1}(D)$ is not reduced then
$f$ is locally a square times a unit. Therefore, for every $x\in X$,
the restriction of $f$ to $T^*_{X,x}$ is either zero or a square
$h_x^2$ where $h_x$ is unique up to a sign. This implies that there is
a ramified cover $\XS\pfeil X$ of degree at most two such that the
pull back of $f$ to $T_\XS^*$ is a square $h^2$. By~\cite{WuM}~6.5, we
have $W_\XS=W_X$. Thus $\CC[\fa^*/W_X]=\CC[\fa^*/W_\XS]$ is integrally
closed in $\CC[T_\XS^*]$. This implies that $f_0$ itself is a square
which is absurd.\qed

\Corollary SmoothCodim1. Let $X$ be affine. Then
$\phi:T^*_X\pfeil\fa^*/W_X$ is smooth in codimension one, i.e., the
set of points where $\phi$ is not smooth is of codimension at least
two.

\Proof: Let $D\subseteq T^*_X$ be a prime divisor. Since $\phi$ is
flat, $D':=\phi(D)$ is a divisor in $S$. Since $D\pfeil D'$ is
generically smooth and $D$ is a reduced component of $\phi^{-1}(D')$
we conclude that $\phi$ is smooth in a generic point of $D$.\qed


\beginsection Appendix A. Appendix A: Lifting and pushing
down differentiable maps

Let $X$ be any smooth manifold and $Y\subseteq X$ a subset. A function
$f$ on $Y$ is differentiable if every $y\in Y$ has an open neighborhood $U$ in
$X$ and a differentiable function $\fQ$ on $U$ such that $f=\fQ|_{Y\cap
U}$. The set of differentiable functions is denoted by $\Cinf(Y)$. We can also
form the sheaf $\Cinf_Y$ of differentiable functions on $Y$. In this paper we
consider only locally closed subsets, i.e., $Y=Z\cap U$ where
$Z\subseteq X$ is closed and $U\subseteq X$ is open. In this case, $U$
is a manifold and an easy argument involving partitions of unity shows
that $\Cinf(Y)$ are just the restrictions of differentiable functions on $U$
to $Y$.

If $K$ is a compact Lie group acting on $X$ let $\pi:X\pfeil X/K$ be
its quotient. Then a function $f$ on $X/K$ is called differentiable if
$f\circ\pi\in\Cinf(X)$. Very often this case can be reduced to the
first one (see ?):

\Theorem. Let $V$ be finite dimensional representation (over $\RR$) of
$K$ and $X\subseteq V$ locally closed, analytic and $K$\_stable. Let
$f_1,\ldots,f_s$ be generators of the ring of $K$\_invariant
polynomials on $V$ and $\pi:V\pfeil\RR^s$ the mapping with components
$f_i$. Then $\pi(X)$ is a locally closed subset of $\RR^s$ and
$X/K\pfeil\pi(X)$ is an homeomorphism and induces an isomorphism
$\Cinf(\pi(Y))\pf\sim\Cinf(X/K)$.

\noindent Note, that Koszul's slice theorem implies that every
$K$\_manifold meets the assumptions of the theorem at least locally.

For any point $x$ of a manifold $M$ let $\hat\cC_x(M)$ be the
completion of $\Cinf(M)$ with respect to the $\fm_x$\_adic
topology. It is a formal power series ring and the image of
$f\in\cC(M)$ in $\hat\cC_x(M)$ is its Taylor series $\hat f_x$. We are
now stating two basic theorems about the relationship of $f$ with its
Taylor series. They are easy consequences of deep theorems by
H{\"o}rmander \cite{Hor}, Malgrange \cite{Mal}, Tougeron \cite{Tou}, and
Bierstone\_Milman \cite{BM}.

Consider the following diagram:
$$
\matrix{&&X\cr&\Links\beta\nearrow&\downarrow\Rechts\pi\cr
U&\pf\alpha&Y\cr}
$$
Given $\alpha$ and $\pi$, we say {\it $\alpha$ lifts to $X$\/} if
there is $\beta$ with $\alpha=\pi\circ\beta$. For every $u\in U$,
$x\in X$ with $y:=\alpha(u)=\pi(y)$ we obtain homomorphisms
$\hat\alpha$ and $\hat\pi$ between completions:
$$
\matrix{&&\hat\cC_x(X)\cr&\Links{\hat\beta}\swarrow&\uparrow\Rechts{\hat\pi}\cr
\hat\cC_u(U)&\pf{\hat\alpha}&\hat\cC_y(Y)\cr}
$$
We say, {\it $\alpha$ lifts formally to $X$\/} if for every $u\in U$
there is $x\in X$ as above and a homomorphism $\hat\beta$ with
$\hat\alpha=\hat\beta\circ\hat\pi$.

Clearly, if $\alpha$ lifts then it lifts formally. A converse is given
by the following theorem:

\Theorem FormLift. We assume
\item\ri $X$ and $Y$ are real\_algebraic manifolds and $\pi$ is a
morphism.
\item\rii There is $h\in\RR[Y]$ such that $\pi$ is a closed embedding
over $Y_0:=\{y\in Y\mid h(y)\ne0\}$.
\item\riii The function $\hq:=h\circ\alpha\in\Cinf(U)$ is non\_zero and
locally analytic.\Par\noindent
Then $\alpha$ lifts if and only if it lifts formally. Moreover, the
lift $\beta$ is unique.

\Proof: The map $\beta$ is unique since by \riii{} the zero\_set of
$\hq$ is nowhere dense and $\pi$ is injective over $Y_0$.

As for existence, since $X$ and $Y$ are algebraic we have to show that
$\RR[Y]\pf{\alpha^*}\Cinf(U)$ extends to
$\RR[X]\pf{\beta^*}\Cinf(U)$.  Because $\alpha$ lifts formally, the
image of $\alpha$ is contained in the image of $\pi$. Thus, if $I$ is
the kernel of $\RR[Y]\pf{\pi^*}\RR[X]$ then
$\alpha^*(I)=0$. Therefore, we obtain a homomorphism
$\alpha^*:R:=\RR[Y]/I\pfeil\Cinf(U)$ and we can think of $R$ as a
subring of $\RR[X]$. Because of \rii, every $f\in\RR[X]$ is of the
form $g/h^N$ where $g\in R$ and $N\in\NN$. Therefore, $\beta^*(f)$ exists
if and only if $\alpha^*(g)$ is divisible by $\hq^N$. Because $\hq$ is
locally analytic, this can be checked by looking at Taylor series
(\cite{Mal}~Thm.~1.1). But for them divisibility holds since $\alpha$
lifts formally.\qed

We consider now the dual situation:
$$
\matrix{X\cr\Links\pi\downarrow&\searrow\Rechts\alpha\cr
Y&\pf\beta&U\cr}
$$
Given $\alpha$ and $\pi$ we say {\it $\alpha$ pushes down to $Y$\/} if
there is a $\beta$ with $\alpha=\beta\circ\pi$. Analogously, {\it
$\alpha$ pushes down formally to $Y$\/} if for every $x\in X$,
$y=\pi(x)$, $u=\alpha(x)$ there is
a homomorphism $\hat\beta$ completing the following diagram:
$$
\matrix{\hat\cC_x(X)\cr\Links{\hat\pi}\uparrow&\nwarrow\Rechts{\hat\alpha}\cr
\hat\cC_y(Y)&{\buildrel\hat\beta\over\leftarrow}&\hat\cC_u(U)\cr}
$$
\Theorem FormPush. We assume
\item\ri $X$ is a manifold and $Y$ is semi\_analytic.
\item\rii $\pi$ is equivariant with respect to an action of a compact
Lie group $K$ and the map on orbit spaces $\pi/K:X/K\pfeil Y/K$ is open.
\item\riii Every orbit $Kx\subseteq X$ has a neighborhood $V$ with an
analytic structure such that $\pi|_V$ is analytic.
\item\riv All fibers of $\pi$ are non\_empty and connected.\Par\noindent
Then $\alpha$ pushes down if and only if it pushes down
formally. Moreover, the push-down $\beta$ is unique.

\Proof: The map $\pi$ is surjective by \riv. So $\beta$ is unique if
it exists. Moreover, $\beta$ exists if and only if $f\circ\alpha$ can
be pushed down for every differentiable function $f$ on $U$. Thus we may
assume $U=\RR$. Since $\pi$ is locally analytic, the set of
smooth points in any fiber of $\pi$ is dense. Since $\alpha$ can be
pushed down formally, the derivatives of its restriction to a fiber of
$\pi$ in a smooth point is zero. Because of \riv, we conclude that
$\alpha$ is constant on the fibers of $\pi$. In particular, $\beta$
exists set\_theoretically.

Let $V\subseteq X$ as in \riii{} which we may choose to be
$K$\_stable. Then $\pi(V)$ is open in $Y$ hence again
semi\_analytic. Let $Z\subseteq Y_0$ be compact. For every $y\in Z$
choose $x\in V$ with $z=\pi(x)$ and a compact $K$\_stable neighborhood
$V_y$ of $Kx$ in $V$. Because of \rii, the image $\pi(V_y)$ is a
neighborhood of $z$. Thus there are finitely many sets $V_i:=V_{y_i}$
such that $Z$ is covered by $\pi(V_i)$. Put
$Z':=\pi^{-1}(Z)\cap(\cup_iV_i)$. Then $Z'$ is a compact subset of $V$
with $\pi(Z')=Z$. This means by definition that $\pi:V\pfeil\pi(V)$ is
semi\_proper.

This means that $\alpha|_V$ satisfies all conditions of
\cite{BM}. Therefore, $\alpha$ can be pushed down to a differentiable function
$\beta_V$ on $\pi(V)$. By uniqueness we have
$\beta_V=\beta_{\pi(V)}$. Because $Y$ is covered by open sets of the
form $\pi(V)$ we have proved that $\beta$ is differentiable.\qed

\beginsection Appendix B. Appendix B: The local structure of
Hamiltonian manifolds

For any $u\in\ft^*_+$ let $L=L(u):=K_u$. A subgroup of this type is
called a Levi subgroup of $K$. Since $L$ contains $T$, its Lie algebra
$\fl$ has a unique complement $\fq$ in $\fk$. Thus we can regard
$\fl^*\subseteq\fk^*$. Let $\fl^r=\{v\in\fl^*\mid K_v\subseteq L\}$
and let $\fl^0$ be its connected component which contains $u$. Both
sets are open in $\fl^*$. Moreover, $K\times^L\fl^0\pfeil\fk^*$ is a
diffeomorphism onto an open subset.

Let $Y$ be a Hamiltonian $L$\_manifold such that
$\mu_L(Y)\subseteq\fl^0$. Then there is a Hamiltonian structure on
$X=K\times^LY$ as follows: let $x=[1,y]\in X$. Then
$T_xX=T_yY\oplus\fq x$ is an orthogonal decomposition and the
symplectic form on $\fq x$ is given by $\omega(\xi x,\eta
x)=\<\mu(y),[\xi,\eta]\>$ for all $\xi,\eta\in\fq$. The moment map on
$X$ is $\mu_X([k,y])=k\mu_Y(y)$.

Now we can state the local cross\_section theorem of Guillemin\_Sternberg:

\Theorem Cross. Let $M$ be a Hamiltonian $K$\_manifold. For a Levi subgroup
$L\subseteq K$ put $M(L):=\mu^{-1}(\fl^0)$ and $M_L:=K\cdot
M(L)$. Assume, $M(L)$ is not empty. Then
\item{\ri} The set $M(L)$ is a Hamiltonian $L$\_manifold: its symplectic form
is the restriction of that of $M$; its moment map is
$M(L)\pfeil\fl^0\into\fl^*$.
\item{\rii} The set $M_L$ is open, dense, connected, and the map
$\Mq:=K\times^L M(L)\pfeil M_L$ is an Hamiltonian isomorphism.
\item{\riii} Let $x\in M$ with $\mu(x)\in\ft^*_+$. Put
$L=K_{\mu(x)}$. Then $Lx$ is an isotropic orbit in $M(L)$.\Par

\noindent
The next theorem describes the neighborhood of an isotropic orbit:

\Theorem Slice. Let $Kx\subseteq M$ be an isotropic orbit. Put $S_x:=(\fk
x)^\perp/\fk x$ which is a symplectic representation of $K_x$. Let
$u_x:=\mu(x)\in(\fk^*)^K$. Then the triple $(K_x,S_x,u_x)$ determines a
neighborhood of $Kx$ uniquely up to Hamiltonian isomorphism.

\noindent Conversely, every such triple occurs: let $H\subseteq K$ be
a closed subgroup, $S$ a symplectic representation of $H$ and
$u_0\in(\fk^*)^K$. We choose a $H$\_invariant unitary structure on
$S$. This induces a moment map $\mu_S:S\pfeil\fh^*$ by
$\mu_S(v)(\xi):={i\over2}\<\xi v,v\>$. Furthermore, choose an
$H$\_stable complement $\fr$ of the Lie algebra $\fh$ in $\fk$ such
that $\fh^*\into\fk^*$. Then $M:=K\times^H(\fh^\perp\oplus S)$ carries
a Hamiltonian structure with moment map
$\mu([k,(u,v)])=u_0+k\cdot(u+\mu_S(v))$. One verifies that the
$K$\_orbit of $x:=[1,(0,0)]$ is isotropic with
$(K_x,S_x,u_x)=(H,S,u_0)$.

\beginrefs


\L|Abk:BM|Sig:BM|Au:Bierstone, E.; Milman, P.|Tit:Composite
differentiable functions|Zs:Ann. of Math. (2)|Bd:116|S:541--558|J:1982||


\B|Abk:GS|Sig:GS1|Au:Guillemin, V.; Sternberg, S.|Tit:Symplectic
techniques in physics|Reihe:-|Verlag:Cambridge University
Press|Ort:Cambridge|J:1984||

\L|Abk:GS2|Sig:GS2|Au:Guillemin, V.; Sternberg, S.|Tit:Multiplicity\_free
spaces|Zs:J. Diff. Geom.|Bd:19|S:31--56|J:1984||

\L|Abk:Ha|Sig:Ha|Au:Hartshorne, R.|Tit:Complete intersections and
connectedness|Zs:Amer. J. Math.|Bd:84|S:497--508|J:1962||

\L|Abk:HeHu|Sig:HH|Au:Heinzner, P.; Huckleberry, A.|Tit:K{\"a}hlerian
potentials and convexity properties of the moment map|%
Zs:Invent. Math.|Bd:126|S:65--84|J:1996||

\L|Abk:HNP|Sig:HNP|Au:Hilgert, J.; Neeb, K.-H.; Plank, W.|%
Tit:Symplectic convexity theorems and coadjoint orbits|%
Zs:Compos. Math.|Bd:94|S:129--180|J:1994||

\L|Abk:Hor|Sig:H{\"o}r|Au:H{\"o}rmander, L.|Tit:On the division of
distributions  by polynomials|Zs:Ark. Mat.|Bd:3|S:555--568|J:1958||

\L|Abk:Ka|Sig:Ka|Au:Karshon, Y.|Tit:Hamiltonian actions of Lie groups|%
Zs:PhD thesis, Harvard|Bd:-|S:-|J:1993||

\L|Abk:KL|Sig:KL|Au:Karshon, Y.; Lerman, E.|Tit:The centralizer of
invariant functions and division properties of the moment map|%
Zs:{\rm To appear in:} Illinois J. Math. {\rm (see {\tt dg-ga/9506008})}%
|Bd:-|S:29 pages|J:1995||

\L|Abk:WuM|Sig:K1|Au:Knop, F.|Tit:Weylgruppe und Momentabbildung|%
Zs:Invent. Math.|Bd:99|S:1-23|J:1990||

\L|Abk:IB|Sig:K2|Au:Knop, F.|Tit:{\"U}ber Bewertungen, welche unter einer
reduktiven Gruppe invariant sind|Zs:Math. Ann.|Bd:295|S:333--363|J:1993||

\L|Abk:CM|Sig:K3|Au:Knop, F.|Tit:The asymptotic behavior of invariant
collective motion|Zs:Invent. math.|Bd:116|S:309--328|J:1994||

\L|Abk:HC|Sig:K4|Au:Knop, F.|Tit:A Harish\_Chandra homomorphism for
reductive group actions|Zs:Annals Math. (2)|Bd:140|S:253--288|J:1994||

\L|Abk:ARC|Sig:K5|Au:Knop, F.|Tit:Automorphisms, root systems, and
compactifications of homogeneous varieties|%
Zs:J. Amer. Math. Soc.|Bd:9|S:153--174|J:1996||

\L|Abk:KP|Sig:KP|Au:Kraft, H.; Procesi, C.|Tit:Closures of conjugacy
classes of matrices are normal|Zs:Invent. Math.|Bd:53|S:227--247|J:1979||

\L|Abk:Le0|Sig:L1|Au:Lerman, E.|Tit:On the centralizer of invariant
functions on a Hamiltonian $G$\_space|%
Zs:J. Differential Geom.|Bd:30|S:805--815|J:1989||

\L|Abk:Le|Sig:L2|Au:Lerman, E.|Tit:Symplectic cuts|Zs:Math. Research
Lett.|Bd:2|S:247--258|J:1995||

\B|Abk:Mal|Sig:Mal|Au:Malgrange, B.|Tit:Ideals of differentiable functions|%
Reihe:Tata Institute of Fundamental Research Studies in Mathematics, No. 3|%
Verlag:Oxford University Press|Ort:London|J:1967||


\L|Abk:Sja|Sig:Sj|Au:Sjamaar, R.|Tit:Convexity properties of the moment map
re-examined|Zs:Preprint {\tt dg-ga/9408001}|Bd:-|S:35 pages|J:1994||

\Pr|Abk:Spr|Sig:Sp|Au:Springer, T.|Artikel:The unipotent variety of a
semisimple group|Titel:Algebraic Geometry (Internat.
Colloq., Tata Inst. Fund. Res., Bombay, 1968)|Hgr:S. Abhyankar, ed.|%
Reihe:-|Bd:-|Verlag:Oxford Univ. Press|Ort:London|S:373--391|J:1969||

\L|Abk:Tou|Sig:Tou|Au:Tougeron, J.-C.|Tit:Fonctions compos{\'e}es
diff{\'e}rentiables: cas alg{\'e}brique|Zs:Ann. Inst. Fourier (Grenoble)|%
Bd:30|S:51--74|J:1980||

\B|Abk:ZaS|Sig:ZS|Au:Zariski, O.; Samuel, P.|Tit:Commutative Algebra II
|Reihe:Graduate Texts in Mathematics~{\bf 29}|Verlag:Springer|%
Ort:Heidelberg|J:1960||

\endrefs

\bye